\documentclass[]{article}
\usepackage{float}
\usepackage{varwidth}
\usepackage{amsfonts}
\usepackage{amsmath}
\usepackage{enumerate}
\usepackage{bm} 
\usepackage{graphicx}
\usepackage{color}
\usepackage{subfig}
\usepackage{booktabs}
\usepackage{lipsum}
\usepackage{microtype}
\usepackage{url}
\usepackage[final]{pdfpages}
\usepackage{breakcites}
\usepackage{array}
\usepackage[english]{babel}
\usepackage{float}
\def\vv{\mathbf v}

\def\dd{\mathbf d}
\def\uu{\mathbf u}

\def\hh{\mathbf h}
\def\Rg{\mathbf R \mathbf g_D}
\def\Rgt{\mathbf R \widetilde{\mathbf g}_D}

\def\UU{\mathbf U}
\def\UUo{\widetilde{\mathbf U}}

\def\DD{\mathbf D}
\def\XX{\mathbf X}

\def\VV{\mathbf V}

\def\Mu{\text{\boldmath$\mu$}}

\def\0a{a_{0}}
\def\t0a{\widetilde{a}_{0}}
\def\1a{a_{1}}

\def\LO{L^2(\Omega)}

\def\review1{\color{red}}

\title{Reduced order models for fluid-structure interaction problems with applications in haemodynamics} 
\author{ Claudia Maria Colciago \footnotemark[4]
\and Simone Deparis\footnotemark[3] % author 4}
}

\begin{document}

\maketitle

\renewcommand{\thefootnote}{\fnsymbol{footnote}}

\footnotetext[3]{CMCS--MATHICSE--SB, \'Ecole Polytechnique F\'ed\'erale de Lausanne, Av. Picard, Station 8, Lausanne, CH--1015, Switzerland.
}

%\footnotetext[6]{CMCS--MATHICSE--SB, \'Ecole Polytechnique F\'ed\'erale de Lausanne, Av. Picard, Station 8, Lausanne, CH--1015, Switzerland}
% \email{luca.dede@epfl.ch}}

%\footnotetext[5]{CMCS--MATHICSE--SB, \'Ecole Polytechnique F\'ed\'erale de Lausanne, Av. Picard, Station 8, Lausanne, CH--1015, Switzerland}
% \email{alfio.quarteroni@epfl.ch}}

\renewcommand{\thefootnote}{\arabic{footnote}}

%\slugger{sinum}{xxxx}{xx}{x}{x--x}%slugger should be set to mms, siap, sicomp, sicon, sidma, sima, simax, sinum, siopt, sisc, or sirev

% history
%\begin{history}
%\received{(Day Month Year)}
%\revised{(Day Month Year)}
%\accepted{(Day Month Year)}
%\comby{(xxxxxxxxxx)}
%\end{history}

% abst
%\begin{abstract}
\section*{Abstract}
This paper deals with fast simulations of the haemodynamics in large arteries
by considering a reduced model of the associated fluid-structure interaction problem, which in turn allows an additional reduction in terms of the numerical discretisation.
The resulting method is both accurate and computationally \emph{cheap}. This goal is achieved by means of two levels of reduction: first, we describe the model equations with a \emph{reduced} mathematical formulation which allows to write the fluid-structure interaction problem as a Navier-Stokes system with non-standard boundary conditions; second, we employ numerical reduction techniques to further and drastically lower the computational costs. The numerical reduction is obtained coupling two well-known techniques: the proper orthogonal decomposition and the reduced basis method, in particular the greedy algorithm. We start by reducing the numerical dimension of the problem at hand with a proper orthogonal decomposition and we measure the system energy with specific norms; this allows to take into account the different orders of magnitude of the state variables, the velocity and the pressure. Then, we introduce a strategy based on a greedy procedure which aims at enriching the reduced discretization space with low \emph{offline} computational costs. 
As application, we consider a realistic haemodynamics problem with a perturbation in the boundary conditions and we show the good performances of the reduction techniques presented in the paper. The gains obtained in term of CPU time are of three orders of magnitude.

\paragraph{Keywords:} Fluid-structure interaction, Navier-Stokes equations, reduced order modelling, proper orthogonal decomposition, reduced basis method, haemodynamics.

%\end{abstract}
%\begin{AMS} \end{AMS}

\pagestyle{myheadings}
\thispagestyle{plain}
\markboth{C.~M.~Colciago,  S.~Deparis, and A.~Quarteroni}{ROM for RFSI}
% main
% 1
\section{Introduction}
When modelling haemodynamics phenomena in big arteries, the resulting model is a complex unsteady fluid-dynamics system, usually coupled with a structural model for the vessel wall. In specific cases, suitable assumptions can be made to reduce the complexity of the model equations. In particular, when the displacement is small, the moving domain can be linearized around a reference steady configuration and the dynamics of the vessel motion can be embedded in the equations for the blood flow. In such way we obtain a \emph{Reduced} Fluid-Structure Interaction (RFSI) formulation where a Navier-Stokes system in a \emph{fixed} fluid domain is supplemented by a Robin boundary condition that represents a surrogate of the structure model.

Although the RFSI model is faster with respect to fully three-dimensional (3D) models where the structure is solved separately, the numerical computation of one heartbeat is still expensive: the resolution of an entire heartbeat, that typically lasts one \emph{physical} second, takes orders of hours of computational time on a supercomputer. A big challenge in realistic applications is to achieve a real time resolution of fluid-structure interaction problems. In particular, in haemodynamics applications, this would grant the possibility to perform real time diagnosis. Nevertheless, the great variability of patient-specific data requires the parametrization of the model with respect to many physical and geometrical quantities. Moreover, as we have recalled above, the complexity of the haemodynamics phenomena requires a mathematical description with complex unsteady models that are difficult to be solved in real time. The RFSI model is already a simpler version of the fully 3D FSI system; a further reduction of the physical model would result in an inaccurate estimation of specific outputs, like e.g. the wall shear stress when using a rigid wall model \cite{Colciago:2012aa,Bazilevs:2007aa}. Thus, to further reduce numerical costs, in this work, we focus on the reduction of the discretization space. In realistic applications, the finite element space has order of $10^6$ degrees of freedom. The aim is to construct a discretization space such that the number of degrees of freedom is reduced to less than $100$ and then to be able to solve one heartbeat in  1 second. 

In the past few years, due to their relevance in realistic applications, a lot of interest has been devoted to discretization reduction techniques for parametrized Partial Differential Equation (PDE) problems (e.g. \cite{Manzoni:2014aa,Lassila:2014aa,Rozza:2008aa,CECAM:2013aa}). These techniques aim to define a suitable \emph{reduced order model} which can be solved with marginal computational costs for different values of the model parameters. Reduced order models are then important in the \emph{many query} context, when a parametrized model has to be solved for different values of the parameter, and in the real time context, when the solution has to be computed with marginal computational costs. To obtain a suitable reduced order model, we typically start from a problem written in a high-fidelity approximation framework, e.g. using the finite element method. The dimension of the discretized system is then drastically reduced through suitable projection operators. The construction of these projection operators is the core of the reduced order technique. Another key concept in the reduction framework is the subdivision of the computational costs into two stages: an \emph{offline} stage, expensive but performed once, and an \emph{online} stage, real time and performed each time new values of the model parameters are considered. During the offline stage the projection space is generated by a reduced basis of functions of the high-fidelity approximation space. 
%The projection matrix defines indeed the \emph{reduced basis} of the new discretization space where a real time resolution can be achieved.

%As an example of reduced order models applied to parametrized PDEs, we cite here the applications to linear elasticity problems when considering variable material properties (see e.g. \cite{Rozza:2008aa, Huynh:2007aa}). %Further applications to convection-diffusion problems can be found in e.g. \cite{Eftang:2011aa, Haasdonk:2008aa}. 
Reduced order models applied to the Burgers equation parametrized with respect to the P\'eclet number is considered in \cite{Yano:2012aa, Nguyen:2009aa}. Other applications of reduced basis techniques applied to fluid problems can be found e.g. in \cite{Wirtz:2014aa, Deparis:2012ab, Iapichino:2012aa, Amsallem:2011aa, Gerner:2011aa, Deparis:2009aa, Grinberg:2009aa, Deparis:2008aa, Gunzburger:2007aa, Sen:2006aa} and in the recent volume \cite{CECAM:2013aa}. 

The aim of this work is indeed to propose a suitable discretization reduction algorithm that can be applied to a RFSI problem. The work is organized as follows. In Section \ref{sec:model} we present the partial differential equations that we are interested in solving. We propose a possible parametrization of the unsteady equations with respect to temporal varying data and with respect to a perturbation of the boundary data. In Section \ref{sec:reduction} we then present how the standard proper orthogonal decomposition algorithm can be applied to the problem at hand in order to generate a suitable reduced space. Moreover, we propose a way to improve the quality of the reduced approximation based on a greedy procedure. Finally, in Section \ref{sec:application} we apply the reduction algorithms presented to a realistic haemodynamic problem. Conclusions follow.

\section{Model equation}\label{sec:model}
Let us consider a three dimensional domain $\Omega$ whose boundary $\partial \Omega$ is divided into three non intersecting parts such that $\partial \Omega = \Gamma \cup \Gamma_{D} \cup \Gamma_N$. 
$\Gamma_D$ is the Dirichlet boundary, typically the inflow of a vessel, $\Gamma_N$ is the Neumann boundary, typically the ouflows, and $\Gamma$ the fluid-structure interface.
We introduce the Hilbert space $V=H^1(\Omega;\Gamma)\:=\{v\in H^1(\Omega)\: v|_{\Gamma} \in H^1(\Gamma)\}$ and the correspondent vectorial spaces $\mathbf V = [V]^3$ and $\mathbf V = [V]^3$. Moreover, we introduce a suitable couple standard finite element spaces $\mathbf V_h$ and $Q_h$ such that $\mathbf V_h \subset \mathbf V$ and $Q_h\subset L^2(\Omega)$ and they represent a stable coupled of finite element spaces for the Navier-Stokes equations. We set $\XX_h := \mathbf V_h \times Q_h$. We define $[t_0\; T]$ a time interval of interest and we divide it into subintervals $[t_{n}\; t_{{n}+1}]$ for $n =0,..,N_T-1$ such that $t_0 < t_1 < t_2< \ldots < t_{N_T} = T$ and $t_{{n}+1} - t_{n} = \Delta t$; let us define $\mathcal N_T =\{0,1,\ldots, N_T, N_T\}$ the collections of all the temporal indexes $n$. For a generic function $\phi(t)$ we use $\phi^{n} := \phi(t_n)$. Finally, we define the operators $D(\cdot)$ and $D_{\Gamma}(\cdot)$ as follows:
\begin{equation*}
D(\vv) = \dfrac{\nabla \mathbf \vv + (\nabla \mathbf \vv)^T}{2} \quad \text{and} \quad D_{\Gamma}(\vv) = \dfrac{\nabla_{\Gamma} \mathbf \vv + (\nabla_{\Gamma} \mathbf \vv)^T}{2}  \qquad \forall \vv \in \mathbf V,
\end{equation*} 
where $\nabla(\cdot)$ is the standard gradient operator and $\nabla_{\Gamma}(\cdot)$ is the tangential component of the gradient with respect to the surface $\Gamma$. 

The RFSI model as presented in \cite{Colciago:2012aa} is an unsteady Navier-Stokes model set on a fixed domain with generalized Robin boundary conditions (For similar models see e.g. \cite{Moireau:2011aa,Figueroa:2009aa,Nobile:2008aa,Figueroa:2006aa}). Let us introduce the velocity and pressure unknowns $[\uu_h, p_h]$ and the corresponding test functions $[\vv_h,q_h]$. Although the RFSI model lives in a fixed domain, it is necessary to define an auxiliary variable which stands for the displacement of the arterial wall $\dd_{s,h}$. Using a backward Euler finite difference method for the time derivatives, the fully discrete weak formulation of the RFSI problem is written as follows: 

for each $n=0,..,N_T-1$, find $[\uu^{n+1}_h, p^{n+1}_h] \in \mathbf X_h$ such that $\uu^{n+1}_h = \mathbf g_D^{n+1}$ on $\Gamma_D$ and
\begin{equation} 
\begin{aligned}
\displaystyle 
&a_0([\uu^{n+1}_h, p^{n+1}_h],[\vv_h,q_h]) + a_1(\uu^{n}_h, \uu^{n+1}_h,\vv_h) =&\\
&\qquad \qquad \qquad \qquad \qquad \qquad F_0(\vv_h; \hh^{n+1}) +
F_{\uu}(\vv_h;\uu^{n}_h) + F_{\dd_s}(\vv_h; \dd^n_{s,h}) \qquad \forall [\vv_h,q_h] \in \mathbf X_h,&
\end{aligned}
\label{eq:5.24.0}
\end{equation}
where
\begin{equation}
\begin{aligned}
	a_0([\uu^{n+1}_h, p^{n+1}_h],[\vv_h,q_h]) = &\int_{\Omega}\bigg(\rho_f\frac{\uu_h^{n+1}}{\Delta t}\cdot \vv_h+ (2\mu \mathbf D(\uu_h^{n+1}) - p_h^{n+1} I ): \nabla\vv_h + q_h^{n+1}\,\nabla\cdot \uu_h \bigg)d\Omega & \\	
	&+ \int_{\Gamma}\bigg(\frac{h_s\rho_s}{\Delta t}\mathbf u_h^{n+1}\cdot\mathbf v_h  + h_s\Delta t\bm \Pi_{\Gamma}(\mathbf u_h^{n+1}):\nabla_{ \Gamma} \mathbf v_h \bigg)d\Gamma,&  \\
	a_1(\uu^{n}_h, \uu^{n+1}_h, \vv_h) = &\int_{\Omega}\rho_f(\uu_h^n \cdot \nabla) \uu_h^{n+1}\cdot \vv_h d\Omega ,&\\
		F_0(\vv_h; \mathbf h^{n+1}) = &\int_{\Gamma_{N}} \mathbf g_N^{n+1} \cdot \mathbf v_h d\Gamma_{N},\\
	F_{\uu}(\vv_h; p^{n}_h) =& \int_{\Omega}\frac{\rho_f}{\Delta t}\mathbf u_h^{n} \cdot \mathbf v_h d\Omega + \int_{\Gamma} \frac{h_s \rho_s}{\Delta t} \mathbf u_h^n \cdot \mathbf v_h d\Gamma,\\
	F_{\dd_s}(\vv_h; \dd^n_{s,n})= & - \int_{\Gamma} h_s \bm \Pi_{\Gamma}(\mathbf d_{s,h}^{n}):\nabla_{\Gamma} \mathbf v_h d\Gamma,
\end{aligned}
\label{eq:5.24}
\end{equation}
with $\dd_{s,h}^{n+1} = \dd_{s,h}^{n} + \Delta t \uu_h^{n+1}$. The parameters $\rho_f$ and $\rho_s$ represent the density of the fluid and solid, respectively; $h_s$ is the thickness of the solid structure surrounding the fluid domain; $\mu$ is the fluid viscosity. $\bm \Pi_{\Gamma}(\cdot)$ is a differential operator with takes into account the structural stiffness and it is defined as $\bm \Pi_{\Gamma}(\vv) = 2\lambda_s D_{\Gamma}(\vv) + \mu_s \nabla_{\Gamma} \cdot \vv I_{\Gamma}$ for $\vv \in \mathbf V$, where $\lambda_s$ and $\mu_s$ are the Lam\'e structural constants and $I_{\Gamma}$  is the identity matrix projected onto the tangential space with respect to the surface $\Gamma$. The functions $g_N(\mathbf x,t)$ and $g_D(\mathbf x,t)$ are sufficiently regular functions that stand for the Dirichlet and Neumann boundary data, respectively. Finally the problem should be equipped with suitable initial condition that, without any loss of generality, we suppose to be equal to zero. 

{\color{black} As said before, the RFSI problem \eqref{eq:5.24.0} is indeed a linearized Navier-Stokes on a fixed domain with a non standard boundary condition on the interface $\Gamma$. In particular is a generalized Robin boundary condition that contains both a mass and a stiffness boundary term to mimic the presence of a compliant arterial wall surrounding the fluid domain (see \cite{Takahito:2014aa} for more detailed on the analysis of partial differential equations with generalized Robin boundary condition). We remark that $\dd_{s,h}^{n}$ does not represents a problem unknown since it is indeed reconstructed as a weighted sum of the velocities at different time instants}
\subsection{Boundary condition}
Problem \eqref{eq:5.24.0} is endowed with Dirichlet velocity boundary condition on the inlet surface $\Gamma_{D}$
. Given the inlet velocity data $\mathbf g_D(\mathbf x,t)$, at the time instant $t_{n+1}$ we impose:
\begin{equation}
\uu_h^{n+1} = \mathbf g_D^{n+1} \qquad \text{ on } \Gamma_{D}.
\label{eq:5.25}
\end{equation}
The Neumann boundary condition $\mathbf D(\uu_h^{n+1})\mathbf n = \mathbf g_N$ is imposed weakly on $\Gamma_N$.
The solution of problem \eqref{eq:5.24.0}-\eqref{eq:5.25} depends on the time variable $t$ through the inlet and outlet data: $\mathbf g_N(\mathbf x,t)$ and $\mathbf g_D(\mathbf x, t)$. We suppose that
\begin{equation}
\mathbf g_D(\mathbf x,t) = \sigma_1(t)\widetilde{\mathbf g}_D(\mathbf x) \quad \text{ and } \quad \mathbf g_N(\mathbf x,t) = \sigma_2(t)\widetilde{\mathbf g}_N(\mathbf x),
\label{eq:5.26}
\end{equation}
that is we separate the contribution of the space and temporal variables in the inlet and outlet data. In realistic applications, the separation of variables \eqref{eq:5.26} often derives directly from modelling choices. If at the outlet we prescribe an average normal stress, no spatial variability is involved in the boundary condition data $\mathbf g_N$. At the inlet, the Dirichlet data is imposed by means of a velocity profile; typically Poiseuille or Womersley profiles are chosen in haemodynamics applications \cite{Formaggia:2009aa}. The separation of variables in $\mathbf g_D(\mathbf x,t)$ in this case is straightforward. 

Assumption \eqref{eq:5.26} allows to write an affine decomposition of the operators in problem \eqref{eq:5.24} with respect to $\sigma_1(t)$ and $\sigma_2(t)$. With respect to the latter we have:
\begin{equation*}
F_0([\vv_h,q_h]; \sigma_2^{n+1},\mathbf g_N) = \sigma_2^{n+1}\int_{\Gamma_{N}} \widetilde{\mathbf g}_N \cdot \mathbf v_h d\Gamma_{out}.
\end{equation*}
The non homogeneous Dirichlet boundary condition \eqref{eq:5.25} is not directly included in the variational form \eqref{eq:5.24.0}. In order to write the affine decomposition with respect to the parameter $\sigma_1(t)$, a suitable choice to embed condition \eqref{eq:5.25} into the variational formulation has to be made. In the literature two possible approaches are proposed: a \emph{strong} imposition, using a lifting function or suitable Lagrange multipliers \cite{Gunzburger:2007aa}, and a \emph{weak} imposition adding suitable penalty variational terms \cite{Sirisup:2005aa, Bazilevs:2007aa}. Due to the fact that the Dirichlet data can be written in the form \eqref{eq:5.26}, a single time independent lifting function can be constructed and properly weighted by a scalar in order to represent the lifting at each temporal instant.

We explain how problem \eqref{eq:5.24.0} is modified when a lifting function for the Dirichlet condition \eqref{eq:5.25} is introduced. Let us directly consider the fully discretized formulation  \eqref{eq:5.24.0}. We define the time independent lifting function $\mathbf R\widetilde{\mathbf g}: \mathbb R^3 \mapsto \mathbb R^3$ such that $\mathbf R\widetilde{\mathbf g} \in \VV_h$ and
\begin{equation*}
\mathbf R\widetilde{\mathbf g} (\mathbf x) = \widetilde{\mathbf g}_D(\mathbf x) \qquad \text{ on } \Gamma_{D}\qquad \text{ and } \qquad\mathbf R\widetilde{\mathbf g}_D (\mathbf x) = 0 \qquad \text{ on } \partial \Omega \backslash \Gamma_{D}.
\end{equation*}
At the time level $t_{n+1}$, the lifting function of the data $\mathbf g_D^{n+1}=\sigma_1^{n+1}\widetilde{\mathbf g}_D$ reads $\mathbf R\mathbf g_D^{n+1} = \sigma_1^{n+1}\mathbf R\widetilde{\mathbf g}_D$.
Then, for each $t_{n+1}$, we introduce the following change of variable:
\begin{equation}
\widetilde{\uu}_h^{n+1} = \uu_h^{n+1} - \Rg^{n+1}.
\label{eq:5.27.1}
\end{equation}
We define the space $\mathbf X_{h,\Gamma_{D}}$ as $\mathbf X_{h,\Gamma_{D}}:= \VV_h\cap [H^1_{\Gamma_D}(\Omega)]^d \times Q_h$ and we observe that $\dd_{s,h}^{n+1} = \sum_{s=0}^{n+1} \Delta t \uu_h^{s} = \sum_{s=0}^{n+1} \Delta t \widetilde{\uu}_h^{s}$ on $\Gamma$.

\subsubsection{Affine decomposition}
Using the definitions of the functionals as in \eqref{eq:5.24}, we are now ready to write the affine decomposition of problem \eqref{eq:5.24.0} with respect to the temporal parameters $\sigma_1(t)$ and $\sigma_2(t)$.
We remark that the lifting function $\Rgt$ does not depend on the time variable, thus the problem parameter at a fixed time level can be gathered in the following vector: 
\begin{equation}
(\Mu^{n+1})^T := [\mu_0,\mu_1,\mu_2] := [\sigma_1^{n+1}, \sigma_2^{n+1}, \sigma_1^{n} ].
\label{eq:5.38}
\end{equation}
{\color{black} One single time step of finite element approximation of the RFSI problem can be written under the form:

for each $n=0,..,N_T-1$, find $\UUo_h^{n+1} \in \mathbf X_{h,\Gamma_{D}}$ such that
\begin{equation} 
\displaystyle 
a(\UUo_h^{n+1}, \mathbf W_h; \UUo_h^{n},\Mu^{n+1} ) = F(\mathbf W_h; \UUo_h^{n}, \DD_h^{n}, \Mu^{n+1})\qquad \forall\mathbf W_h\in \mathbf X_{h,\Gamma_{D}},
\label{eq:5.35}
\end{equation}
where 
\begin{equation}
\begin{aligned}
&a(\UUo_h^{n+1}, \mathbf W_h; \UUo_h^{n}, \Mu^{n+1} ) :=\0a(\UUo_h^{n+1}, \mathbf W_h ) + \mu_2 a_1( \Rgt, \UUo_h^{n+1},\mathbf W_h)+ a_1(\UUo_h^n, \UUo_h^{n+1},\mathbf W_h),&\\
&F(\mathbf W_h; \UUo_h^{n},\dd_{s,h}^{n}, \Mu^{n+1})):=\mu_1 F_0(\mathbf W_h;\widetilde \hh) + F_{\uu}(\mathbf W_h; \UUo_h^{n}) + \mu_2 F_{\uu}(\mathbf W_h; \Rgt) + F_{\dd_s}(\mathbf W_h; \dd_{s,h}^{n})&\\
&\qquad \qquad \qquad- \mu_0 a_0(\Rgt,\mathbf W_h) - \mu_0 a_1( \UUo_h^{n}, \Rgt,\mathbf W_h) - \mu_2 \mu_0 a_1( \Rgt, \Rgt,\mathbf W_h).&
\end{aligned}
\label{eq:AffDec}
\end{equation} 
Due to the fact that we use a semi-implicit treatment of the convective term the formulation of the RFSI problem at one single time instant $t_{n+1}$ can be interpreted as a linear steady problem parametrized with respect to $\Mu^{n+1}, \UUo_h^n$ and $\dd_{s,h}^n$.

Furthermore, we can introduce a parameter in the inlet flow rate function representing a small perturbation with respect to a reference value: the inlet flow rate function \eqref{eq:5.26} is then modified as 
\begin{equation}
\mathbf g(\mathbf x,t; \alpha) = \theta(\alpha, t) \sigma_1(t)\widetilde{\mathbf g}_D(\mathbf x),
\label{eq:PerturbG}
\end{equation}
where $\alpha \in \mathcal D$, being $\mathcal D$ the set of the admissible value of $\alpha$. {\color{black} The same affine decomposition \ref{eq:AffDec}, with the following modification: the parameter becomes 
$ (\Mu^{n+1})^T := [ \mu_0, \mu_1,\mu_2, \mu_3] := [ \sigma_1^{n+1}, \sigma_2^{n+1}, \sigma_1^{n}, \theta^n(\alpha) ]$ and in \eqref{eq:AffDec} we substitute $\mu_0$ with $\mu_3\mu_0$ and $\mu_2$ with $\mu_2\mu_3$.}

\section{Numerical reduction}\label{sec:reduction}
In this section we briefly introduce some of the basic concepts of the reduced basis method that are useful to our purpose. For more details on the reduced basis theory we address the interested reader to e.g. \cite{Rozza:2008aa,hesthaven2016certified,quarteroni2015reduced}. We already introduced $\UUo_h^{n+1}$ that, at each time instant is the a high-fidelity approximation of the exact solution and is computed as a finite element solution with a sufficiently fine mesh. The solutions $\UUo_h^{n+1}$ of problem \eqref{eq:5.35} are, in general, expensive to obtain from the computational point of view, since in realistic applications the finite element spaces has order of  $10^6$ degrees of freedom and the complexity of the geometrical domain does not always allow for the generation of structural meshes. We conclude that due to the magnitude of the finite element problem a real time computation would be impossible to achieve. 

As in the standard reduced basis theory, we state the following assumption: the family of solutions $\UUo_h^{n+1} = \UUo_h^{n+1}(\Mu^{n+1})$ obtained for different realizations of the parameters belongs to a low dimensional manifold $\mathcal M_{h}^{\Mu}$. The aim of the reduction techniques is to find a suitable approximation of the manifold $\mathcal M_h^{\Mu}$ through the construction of a low dimensional space $\XX_{N}\subset \XX_{h,\Gamma_D}$. The dimension of the reduced space $N$ needs to be orders of magnitude lower that the dimension of the finite element space. The reduced approximation of RFSI problem reads: 

given $\UUo_N^{0} = \UUo_h^{0}$, for each $n=0,..,N_T-1$, find $\UUo_N^{n+1} \in \mathbf X_{N}$ such that
\begin{equation} 
\displaystyle 
a(\UUo_N^{n+1}, \mathbf W_N; \UUo_N^{n}, \Mu^{n+1} ) = F(\mathbf W_N; \UUo_N^{n}, \DD_N^{n}, \Mu^{n+1}) \qquad \forall\mathbf W_N\in \mathbf X_N,
\label{eq:5.36}
\end{equation}
where $a(\cdot,\cdot)$ and $F(\cdot)$ are defined as in \eqref{eq:AffDec}. 

\subsection{Proper Orthogonal Decomposition}\label{sec:POD}
We apply a discretization reduction to the RFSI problem \eqref{eq:5.35} and the Proper Orthogonal Decomposition (POD) method. In the context of this work we only detail the specific choices performed in relation to the problem at hand, for more details about POD applied to fluid problems we address the reader to e.g. \cite{Rowley:2005aa, Willcox:2002aa}. 

We define a subset of temporal indexes $\mathcal N_S \subset \mathcal N_T$ with cardinality $N_S$ and consider the solutions of problem \eqref{eq:5.24.0} at the time instants $t^{n_S}$ for $n_S \in \mathcal N_S$. The solutions $\UUo^{n_S}_h$, called \emph{snapshots}, represent our starting point for the POD analysis. Since the RFSI problem \eqref{eq:5.35} is a saddle point problem in two variables (velocity and pressure) with different characteristic order of magnitude, we split the POD into two eigenvalue decompositions: one for the velocity variable and another for the pressure one \cite{Gerner:2012aa}. We measure the energy associated to the snapshots using the following scalar products: for the velocity, we set
\begin{equation}
(\uu_h,\vv_h)_{\VV} := (\uu_h,\vv_h)_{\mathbf H^1(\Omega)} + (\uu_h,\vv_h)_{\mathbf H^1(\Gamma)}, \qquad \forall \uu_h, \vv_h \in \VV_h \subset \VV ( = \mathbf H^1_{\Gamma_D}(\Omega;\Gamma) ),
\end{equation} 
and for the pressure,
\begin{equation}
(p_h,q_h)_{Q} := (p_h,q_h)_{L^2(\Omega)}, \qquad \forall p_h,q_h \in Q_h\subset Q:=L^2(\Omega).
\end{equation}
Then, we compute the two Gramian matrices
\begin{equation}
G_{ij}^{\uu} = (\uu_h^i,\uu_h^j )_{\VV} \quad \text{ and }\quad G_{ij}^{p} = (p_h^i,p_h^j )_{Q}  \qquad \forall j,i \in \mathcal N_S,
\end{equation}
and we perform the eigenvalue decomposition of $G^{\uu}$ and the one of $G^{p}$, obtaining the pairs $(\lambda_k^{\uu}, \bm \zeta_k^{\uu})$ and $(\lambda_k^{p}, \bm \zeta_k^{p})$ where $\lambda_k^{\uu}, \lambda_k^{p} \in \mathbb R$ and $\bm \zeta_k^{\uu}, \bm \zeta_k^{p} \in \mathbb R^{N_S}$ are the $k-th$ eigenvalues and eigenvectors of the velocity and pressure Gramian matrices, respectively, for $k \in \mathcal N_S$. Fixing the same tolerance for both the velocity and pressure decompositions, we select the first $N^{\uu}$ and $N^{p}$ eigenpairs such that:
\begin{equation}
\dfrac{\sum_{j=1}^{N^{\uu}}\lambda^{\uu}_j}{\sum_{k=1}^{N_S}\lambda_k^{\uu}} \geq 1 - tol \quad \text{ and } \quad \dfrac{\sum_{j=1}^{N^p}\lambda^{p}_j}{\sum_{k=1}^{N_S}\lambda_k^p} \geq 1 - tol,
\label{eq:5.33}
\end{equation}
respectively.
The $j-$th velocity eigenfunction $\bm \phi^{\uu}_j \in \VV_h$ is reconstructed using the linear combination:
\begin{equation*}
\bm \phi_j^{\uu} = \frac{1}{\lambda_j} \sum_{n_S \in \mathcal N_S} [\bm \zeta^{\uu}_j]_{n_S} \uu_h^{n_S}, \quad \text{ for } j=1,..,N^{\uu}.
\end{equation*}
Similarly for $\phi^{p}_j \in Q_h$ for $j=1,..,N^{p}$. {\color{black} We remark that, since the velocity basis are linear combinations of solutions of problem \eqref{eq:5.35}, they all verify $\int_{\Omega} q_h \nabla \cdot \bm \psi^{\uu}_{j} =0$, $\forall q_h \in Q_h$ for $j=1,..,N^{\uu}$. Thus, the linear system induced by the bilinear form $a(\cdot,\cdot)$ as in \eqref{eq:5.24} would be singular if we consider the functional spaces generated from the velocity functions $\bm \phi^{\uu}_{j}$ and the pressure modes $\phi^{p}_j$. One of the possibility often employed in the context of Navier-Stokes equations is to restrict the system and to solve the problem only for the velocity unknown (see e.g. \cite{Burkardt:2006aa}). Unfortunately, this is not possible when considering problem \eqref{eq:5.24}. The generalized boundary condition applied on $\Gamma$ derives from a structural model which solution is driven by the pressure condition set on the external boundary in the structural model (see \cite{Colciago:2012aa, Moireau:2011aa}). If we solve the reduced system not taking into account the pressure variable, we cannot recover the velocity on the boundary $\Gamma$ and the output functionals that depends on these values (e.g wall shear stress).} For these reasons, following \cite{Rozza:2007aa}, for each selected pressure mode $\phi_j^{p}$, we define the corresponding supremizer function $\bm \sigma_j \in \VV_h$ as the solution of the following problem:
\begin{equation}
(\bm \sigma_{j}, \vv_h) = \int_{\Omega} \phi_j^{p} \nabla \cdot \vv_h d\Omega \qquad \forall \vv_h \in \VV_h, \qquad \text{ for } j=1,..,N^p.
\label{eq:5.34}
\end{equation}
We then add them to the POD basis. The POD reduced space $\XX_N^{POD}$ associated to the RFSI model is  composed by  the basis functions $\{\bm \psi_j\}_{j=1}^{N^\uu + 2 \times N^p}$, $\bm \xi_j \in \XX_h$ defined as follows:
\begin{equation}
\begin{aligned}
&\bm \psi_j= [\bm \phi_j^{\uu},\, 0]^T \qquad & & \text{ for } j= 1,..,N^{\uu} \\
&\bm \psi_{N^{\uu}+j} = [\mathbf 0,\,\phi_j^{p}]^T \qquad & & \text{ for } j= 1,..,N^p \quad \text{and}  \\
& \bm \psi_{N^p + N^{\uu} +j}= [\bm \phi^{\bm \sigma}_j,\, 0]^T \qquad & & \text{ for } j= 1,..,N^p, 
\end{aligned}
\label{eq:finalBasis}
\end{equation}
where $\bm \phi^{\bm \sigma}_j$ for $j= 1,..,N^p$ represent the orthonormalization of the supremizer functions  $\bm \sigma_j$, obtained with a Gram-Schmidt algorithm with respect to the scalar product $(\cdot,\cdot)_{\VV}$.}

\subsection{Greedy enrichment}\label{sec:Greedy} 
The bottleneck of the POD procedure is the computation of the high-fidelity solutions $\UUo_h^n$ necessary to build the correlation matrix: we have to solve a finite element problem $N_T$ times. Moreover if we choose $\mathcal N_S = \mathcal N_T$, the Gramian matrix becomes too large its eigenvalue decomposition gets too much expensive. We can envision two situations where we would like improve the quality of the approximation obtained with the POD reduced space without changing the snapshots sample. For example, if
$N_S$ is five smaller than $N_T$,  the \emph{information} carried by the snapshots sample refers to only the $25\%$ of the entire set of the truth solutions. Is it possible to improve the quality of the reduced approximation, without increasing the number of snapshots selected? In another scenario, suppose that a perturbation parameter $\alpha$ is introduced in the unsteady problem \eqref{eq:5.35}, as proposed in \eqref{eq:PerturbG}, and that the snapshots are computed for a specific value of $\alpha=\alpha_1$.
We would like to generate a reduced space that suitably approximates also the truth solutions for other values of $\alpha$ without recomputing all the high-fidelity snapshots.

With these two scenarios in mind, we propose a strategy to improve the quality of the reduced approximation based on a greedy algorithm. For references to standard greedy algorithms applied to parametrized PDEs see e.g. \cite{hesthaven2016certified,quarteroni2015reduced}.

%When dealing with realistic applications, we have to take into account also the memory costs of storing a large number of truth solutions. 
%Let us suppose to use the POD algorithm to reduce the temporal variability for an evolution problem. Typically, when solving an unsteady problem for all the time levels $t_n = n\Delta t$ with $n=0,..,N_T$, we do not store all the solutions $u^n_h$ for $n=0,..,N_T$, but only a subset of them. The number and the distribution of the stored solutions have to be chosen such that they do not affect the post-processing phase, when we aim to evaluate specific output functionals of the physical system. 

We introduce another solution $\UU_{N,h}^n$ that belongs to an \emph{intermediate} problem between \eqref{eq:5.35} and \eqref{eq:5.36}: find $\UU_{N,h}^n \in \XX_h$ such that
\begin{equation}
\displaystyle 
a(\UU_{N,h}^{n+1}, \mathbf W_h; \UUo_N^{n},\Mu^{n+1} ) = F(\mathbf W_h; \UUo_N^{n}, \DD_N^{n}, \Mu^{n+1})\qquad \forall\mathbf W_h\in \mathbf X_{h,\Gamma_D},
\label{eq:5.37}
\end{equation} 
We notice that problem \eqref{eq:5.37} is set in the high-fidelity approximation framework but the right hand side and the advection field are defined by \eqref{eq:5.36}. In fact, in \eqref{eq:5.35}, these terms are evaluated using the truth solution $\UUo_{h}^n$, while in \eqref{eq:5.37} it is evaluated using the reduced solution $\UUo_N^n$, as in problem \eqref{eq:5.36}. Considering the error between $\UUo_N^n$ and $\UUo_h^n$ in a generic norm $\|\cdot\|_*$, the following triangular inequality holds:
\begin{equation*}
\|\UUo_N^{n+1} - \UUo_h^{n+1}\|_{*} = \|\UUo_N^{n+1} - \UU_{N,h}^{n+1} + \UU_{N,h}^{n+1} - \UUo_h^{n+1}\|_{*} \leq \|\UUo_N^{n+1} -\UU_{N,h}^{n+1}\|_{*} + \|\UU_{N,h}^{n+1} - \UUo_h^{n+1}\|_{*}
\end{equation*}
The greedy procedure that we propose focuses on the first contribution $\|\UUo_N^{n+1} -\UU_{N,h}^{n+1}\|_{*}$. 
Subtracting problem \eqref{eq:5.36} from \eqref{eq:5.37} allows to state a result of Galerkin orthogonality:
\begin{equation*}
a(\UU_{N,h}^{n+1} - \UUo_N^{n+1}, \mathbf W_h; \UUo_N^{n},\Mu^{n+1} ) = 0.
\end{equation*}
We assume that the dual norm of the residual can be used as an indicator of the error $\|\UU_N^{n+1} - \UU_h^{n+1}\|_{\XX}$. In particular, at each time level $t_{n+1}$, we consider
\begin{equation}
r_N^{n+1}(\mathbf W_h):= F(\mathbf W_h; \UUo_N^{n},\DD_N^{n}, \Mu^{n+1}) - a(\UUo_N^{n+1}, \mathbf W_h; \UUo_N^{n}, \Mu^{n+1})
\label{eq:5:32}
\end{equation}
and its associated dual norm $\|r_N^{n+1}(\mathbf W_h)\|_{\XX'}$.

We now have defined all the necessary quantities, we can proceed presenting the steps to be performed when we want to enrich the POD basis with a greedy algorithm. {\color{black}First, perform a POD on the snapshots $\UUo_h^{n_S}$, for $n_S\in\mathcal N_S$ and we construct the reduced space $\XX_N^{POD}$. Then, we start the greedy enrichment setting $\XX_N=\XX_N^{POD}$:
\begin{itemize}
\item[1.] Generate the reduced basis solutions $\UUo_N^n$, $n \in \mathcal N_T$, by solving the reduced order problem \eqref{eq:5.36}.
\item[2.] Compute the dual norms of the residuals $\|r_N^{n}(\mathbf W_h)\|_{\XX'}$, $n \in \mathcal N_T$, which are used as error indicators. 
\item[3.] Select $n^*$ such that
\begin{equation*}
n^* = \arg \max_{n \in \mathcal N_T} \|r_N^{n}(\mathbf W_h)\|_{\XX'}.
\end{equation*}
\item[4.] Compute the $\UU_{N,h}^{n^*}$ by solving the reduced order problem  \eqref{eq:5.37}.
\item[5.] Split $\UU_{N,h}^{n^*}$ into its velocity and pressure components, $\uu_h^{n^*}$ and $p_h^{n^*}$, respectively. Compute the supremizer $\bm \sigma_{n^*}$ associated with the pressure component.
\item [6.] Compute $\bm \phi^{\uu}$ representing the orthonormalization of the velocity function $\uu_h^{n^*}$ with respect to the reduced space $\XX_N$, obtained with a Gram-Schmidt algorithm considering the scalar product $(\cdot,\cdot)_{\VV}$; similarly for $\phi^{p}$ and $p_h^{n^*}$.
\item [7.] Build $\XX_{N+2}=\XX_{N}\oplus \{\bm \psi^{\uu},\bm \psi^{p}\}$ defined as is \eqref{eq:finalBasis}.
\item [8.] Compute $\bm \phi^{\bm \sigma}$ representing the orthonormalization of the velocity function $\bm \sigma^{n^*}$ with respect to the reduced space $\XX_{N+2}$, obtained with a Gram-Schmidt algorithm considering the scalar product $(\cdot,\cdot)_{\VV}$.
%= \dfrac{\uu_h^{n^*} - \Pi_{X_N}\uu_h^{n^*}}{\|\uu_h^{n^*} - \Pi_{X_N}\uu_h^{n^*}\|_{\VV}}$ and similarly for $\phi^{p}$ and $\bm \phi^{\bm \sigma}$
\item [7.] Build $\XX_{N+3}=\XX_{N+2}\oplus \{\bm \psi^{\bm \sigma}\}$ defined as is \eqref{eq:finalBasis}.
\item[8.] Update the structures for the online computation of the reduced solutions and the dual norms of the residuals.
\item[9.] Set $N = N+3$ and $\XX_N = \XX_{N+3}$. Repeat until a predefined stopping criterion is satisfied.
\end{itemize}

\textbf{Remark.} We remark that the functions that are added to the space $X_N$ in step 5 are derived from $\UU_{N,h}^{n^*}$ and not the truth solution $\UUo_h^{n^*}$. We have no guarantee that $\UU_{N,h}^{n^*}$ is close to $\UUo_h^{n^*}$ or that it belongs to the low dimensional manifold $\mathcal M_{h}^{\Mu}$ of the truth solutions. We would like also to remark that even if we are trying to reduce the error  $\UUo_N^{n+1} -\UU_{N,h}^{n+1}$, to date we have no proof that the algorithm converges. In fact, we cannot theoretically prove that
\begin{equation}
\||\UUo_N^{n+1} -  \UUo_h^{n+1}\||_{*} \leq \||\UUo_{N-1}^{n+1} -  \UUo_h^{n+1}\||_{*}. 
\end{equation}
For this lack of theoretical convergence results, to stop the greedy enrichment procedure, we rather opt for a fixed number of solutions $N_{max}$ chosen a priori,
instead of using a certain tolerance on the a posteriori error estimator. Nevertheless, in the next section we will show some numerical evidence that the greedy enrichment is able to improve the quality of the approximation space.

\section{Application to a femoropopliteal bypass}\label{sec:application}
\subsection{Application and motivation}
Atherosclerotic plaques often occur in the femoral arteries. The obstruction of the blood flow results in a lower perfusion of the lower limbs and the most common symptom of this disease is an intermittent claudication, which affects the 4$\%$ of people over the age of 55 years \cite{Giordana:2005aa}. In order to restore the physiological blood circulation, different medical treatments are possible. 
%Depending on the severity of the stenosis, chemical drugs and/or stents can be used to reduce the deposit of fat inside the arterial wall. 
In critical cases, the stenosis is treated with surgical intervention that \emph{bypasses} the obstruction using a \emph{graft} and providing an alternative way where blood can flow. The bypass creates a side-to-end anastomosis between the graft and the upstream artery (\emph{before} the occlusion) and an end-to-side anastomosis with the distal downstream part. 
%The graft can be autologous (a portion of a vein taken from the same patient) or it can be artificially made. Despite the fact that autologous venous grafts represent the principal treatment for serious stenoses in peripheral arteries, the rate of failure of such a procedure continues to be high \cite{Giordana:2005aa}. 
In particular, the design of end-to-side anastomosis affects the flow downstream the bypass and provokes remodelling phenomena inside the arterial wall. The arteries adapt their size in order to maintain a certain level of shear stress, which results in a thickening of the intima layer and in an increasing risk of thrombi formation. The arterial wall remodelling is in fact linked with haemodynamic factors such as the wall shear stress magnitude and direction. Moreover, velocity profiles and separation of flows have been investigated when studying the bypass end-to-side anastomosis \cite{Loth:2002aa, Loth:2008aa}. Studies with idealized geometrical models have been proposed in order to define an optimal design for the anastomosis \cite{Lassila:2013aa}. Nevertheless, the geometry of the vessel is one of the most important factors that affect the pattern of the wall shear stress. Further, patient specific data would be required in order to analyse each particular case.

\begin{figure}[t]
    \begin{center}
    \subfloat[End-to-side geometry]{
      \includegraphics[width=4cm]{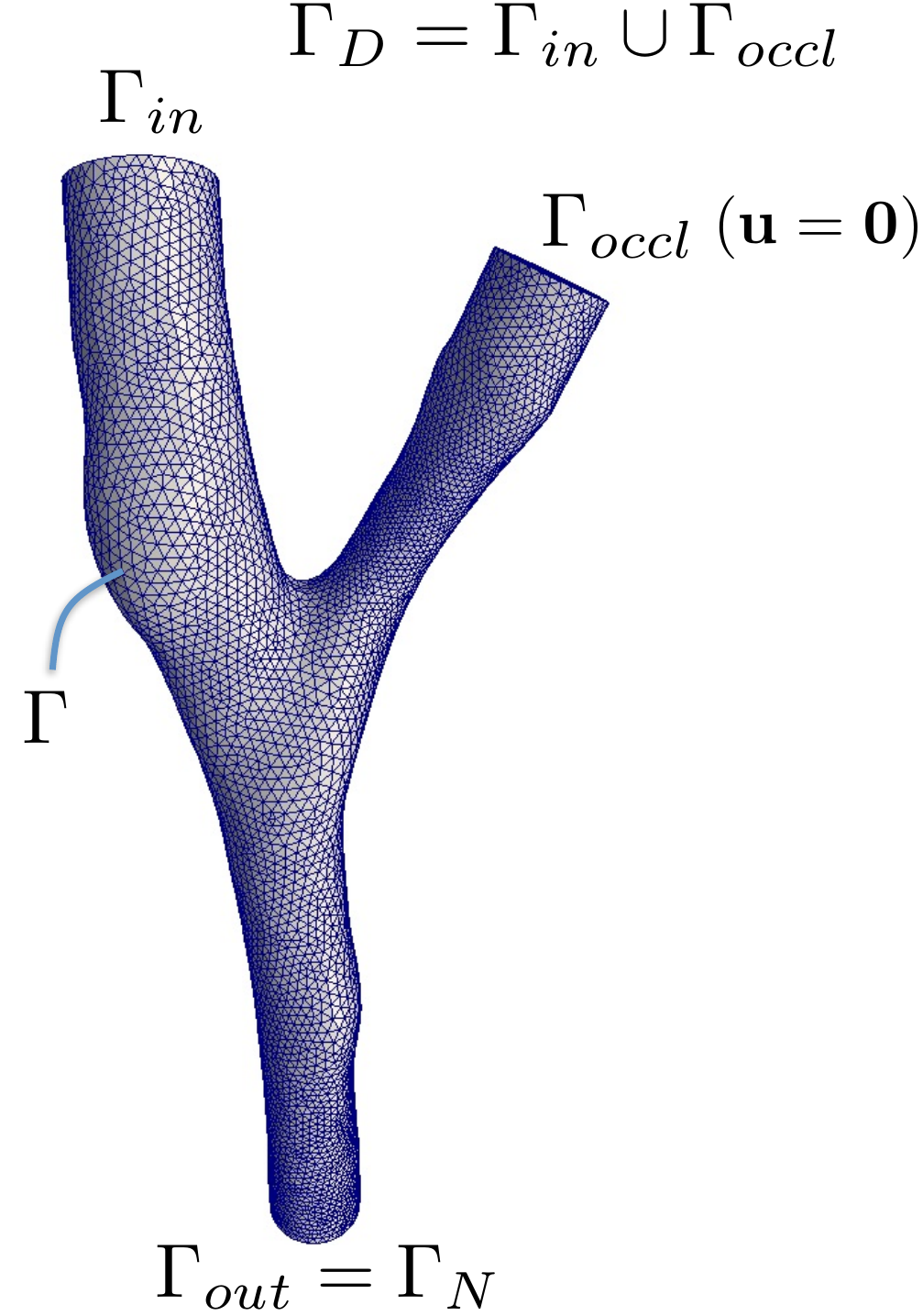}}
      \subfloat[Inlet flow rate]{
      \includegraphics[width=5.3cm]{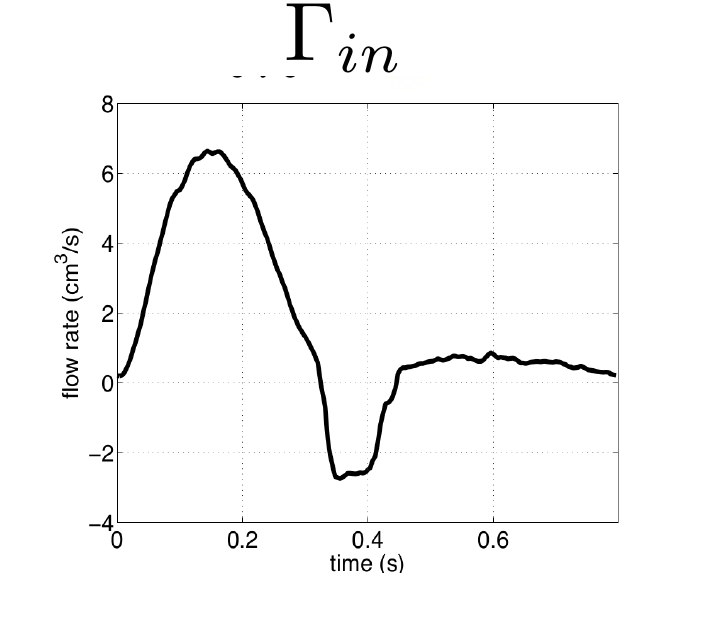}}
      \subfloat[Outlet pressure]{
      \includegraphics[width=5cm]{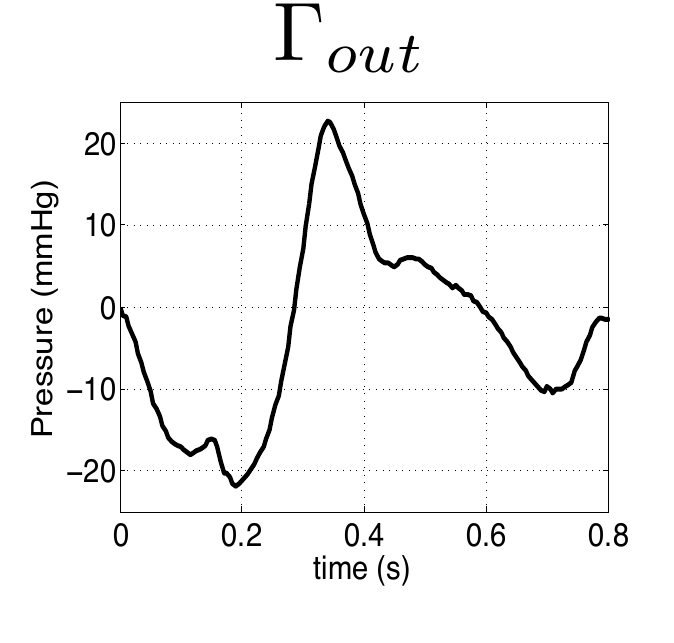}}
    \end{center}
    \caption{Realistic geometry of the end-to-side anastomosis. Graph of the inlet and outlet boundary conditions.}
    \label{fig:BypassAnast}
\end{figure}

We focus our attention on the patient-specific femoropopliteal bypass performed with a venous graft bridging the circulation from the femoral artery to the popliteal one. As a domain of interest we select the end-to-side anastomosis (see Figure \ref{fig:BypassAnast}). The geometry was reconstructed by CT-scan images as it is detailed in \cite{Marchandise:2011aa} and inlet and outlet flow rates are provided from the experimental data in \cite{Marchandise:2011aa}. We compute the Reynolds number as $Re = \dfrac{4\rho_f Q_{in}}{\pi D \mu}$, being $\rho_f$ the blood density, $Q_{in}$ the inlet flow rate, $D$ the vessel diameter and $\mu$ the blood viscosity. The average Reynolds number ranges from 144 and 380 (at the systolic peak), in agreement with the values provided in \cite{Loth:2008aa}.

\subsection{Test case}\label{sec:TestCase}
\subsubsection{Application of the POD algorithm}
In this section we investigate the behaviour of the POD and the greedy enriched POD algorithms on a case representing the femoropopliteal bypass application where the finite element resolution is performed on a coarse mesh. The usage of a coarse grid allows us to lower the offline computational costs and, thus, to test and compare several reduced basis approximations. Since we are interested in the realistic application of the femoropopliteal bypass, the physical parameters and boundary data are patient-specific. The coarse mesh is composed by 5'823 tetrahedra and 1'309 vertices. To obtain the high fidelity solutions of the RFSI model we use standard $\mathbb P^1$+Bubble-$\mathbb P^1$ finite elements for a total of 22'702 degrees of freedom.
{\color{black} The boundary conditions are periodic with period of 0.8 s (one heartbeat). We set the solutions at time t=0 equal to zero. To get rid of the dependence of these initial condition we perform the simulation of an entire heartbeat and we focus on the solutions obtained for the subsequent heartbeat.} Thus, to test the POD reduction algorithm, we compute the high fidelity numerical solutions for a time lapse corresponding to the second heartbeat, from $t_0=0.8$ s to $t_{N_T}=1.6$ s with a time step $\Delta t=0.001$ for a number of time intervals $N_T =800$. We denote with the superscript $n \in \mathcal N_T$ varying from 0 to $N_T$ the sequence of computed solutions:
\begin{equation*}
\UU_h^{n} \approx \UU_h(t_n) \quad \text{ where } t_0 = 0.8, t_1=t_0 + \Delta t, t_2=t_0 + 2 \Delta t, t_3=t_0 + 3 \Delta t, .. , t_{N_T} = 1.6 \text{s}.
\end{equation*}
We save the finite element solutions every five time steps and we use the apex $n_S \in \mathcal N_S$, $n_S= 5k$, with $k=0,..,N_S$ ($N_S = 160$) to address the stored functions, that represent the snapshot sample:
\begin{equation*}
\UU_h^{n_S} \approx \UU_h(t_{n_S}) \quad \text{ where } t_{5} = 0.805, t_{10}=0.810, t_{15}=0.815,t_{20}=0.820, .. , t_{n_{N_S}} = 1.6 \text{s}.
\end{equation*}
Indeed, we compute the POD starting from the 160 snapshots $\UU_h^{n_S}$, $n_S \in \mathcal N_S$, which represent the $25\%$ of the finite element solutions computed for the second heartbeat. To check the quality of the reduced space approximations, we monitor the following errors:
\begin{itemize}
\item  relative error of the velocity at time $t_{n_S}$ and correspondent space-time error:
\begin{equation}
\varepsilon_N(\uu^{n_S}):=\dfrac{\|\uu_N^{n_S} - \uu_h^{n_S}\|_{\VV}}{\|\uu_h^{n_S}\|_{\VV}} \quad \text{and} \quad E_N(\uu):=\dfrac{\bigg(\sum_{n_S \in \mathcal N_S}\big( \|\uu_N^{n_S} - \uu_h^{n_S}\|_{\VV}\big)^2\bigg)^{1/2}}{\bigg(\sum_{n_S\in \mathcal N_S}\big( \|\uu_h^{n_S}\|_{\VV}\big)^2\bigg)^{1/2}};
\end{equation}
\item  relative error of the pressure at time $t_{n_S}$ and correspondent space-time error:
\begin{equation}
\varepsilon_N(p^{n_S}):=\dfrac{\|p_N^{n_S} - p_h^{n_S}\|_{\LO}}{\|p_h^{n_S}\|_{\LO}}\quad \text{and} \quad 
E_N(p):=\dfrac{\bigg(\sum_{n_S \in \mathcal N_S}\big( \|p_N^{n_S} - p_h^{n_S}\|_{\LO}\big)^2\bigg)^{1/2}}{\bigg(\sum_{n_S \in \mathcal N_S}\big( \|p_h^{n_S}\|_{\LO}\big)^2\bigg)^{1/2}};
\end{equation}
\item space-time dual norm of the residual scaled with respect to the space-time norm of the global solution
\begin{equation}
R_N(\UU):=\bigg(\dfrac{N_S}{N_T}\bigg)^{1/2}\dfrac{\bigg(\sum_{n\in \mathcal N_T}\|r_N^{n}(\mathbf W_h)\|_{\XX'}^2\bigg)^{1/2}}{\bigg( \sum_{n_S \in \mathcal N_S}\|\UU_h^{n_S}\|_{\XX}\big)^2\bigg)^{1/2}}; 
\end{equation} 
\end{itemize}

We build a sequence of POD reduced spaces with decreasing values of the tolerance $tol$ and we compute the aforementioned indicators for each one of the reduced spaces generated. The space-time errors are reported in Table \ref{tab:6.1}. In particular, we show: the number of selected velocity modes ($\# \uu$  basis); the number of selected pressure modes ($\# p$  basis); the total number of basis functions composing the reduced space ( $\#$ basis = $\# \uu$  basis + $2 \times \# p$  basis ); the space-time errors and residuals as defined above. Since the problem at hand is unsteady and the solution at a time instant $t_n$ depends on the solutions at the previous instants, the POD model errors $E_N(\uu)$ and $E_N(p)$ are bounded from above by the fixed tolerance but they are however of the same order of magnitude (see Table \ref{tab:6.1}). 
%In Figure \ref{fig:6.2} we display the global absolute errors $\|\UU_N^{n_S} - \UU_h^{n_S}\|_{\XX}$ and the dual norm of the residuals ($\|r_N^{n}(\mathbf W_h)\|_{\XX'}$). 
We notice that, even if $\|r_N^{n}(\mathbf W_h)\|_{\XX'}$ does not represent an upper bound for the error, nevertheless, from experimental results, we can use it as an indicator of $\|\UU_N^{n_S} - \UU_h^{n_S}\|_{\XX}$ (see Figure \ref{fig:6.2}). {\color{black} The apparent strong correlation between the dual norm of the residual and the global error norm is probably due to the strong contribution of the mass term in the unsteady formulation. Indeed, if we choose a time step of 0.001, the mass matrix is multiplied for a factor of $10^3$.} We remark that the magnitudes of the absolute errors for the velocity span from $10^{-1}$ to $10^2$ and the associated velocity solutions norms are of order of $10^{2}-10^{3}$. For the pressure, we have absolute errors of order $10^{0}-10^{3}$, while their solutions norms are of order $10^{3-5}$. {\color{black} Indeed, in absolute terms the global error is mostly related with the pressure one.}

\begin{table}[H]
\centering
\begin{tabular}{cccccccc}
\toprule
	$tol$ & $tol^{1/2}$ & $\# \uu$  basis & $\# p$ basis & $\#$ basis & $E_N(\uu)$ & $E_N(p)$ & $R_N(\UU)$ \\ 
\midrule
	$1e-2$ & $1e-1$ 	    &  8 &  1 &  10  & 2.34e-1   & 3.73e-2   & 4.17e-2\\
	$1e-3$ &$3.16e-2$ 	& 16 & 	1 &  18  & 2.16e-1   & 4.50e-2   & 4.92e-2\\
	$1e-4$ &$1e-2$ 	    & 28 &  2 &  32  & 6.83e-2   & 1.32e-2   & 9.92e-3\\
	$1e-5$ &$3.16e-3$ 	& 44 & 	3 &  50  & 9.09e-3   & 1.89e-3   & 1.42e-3\\
	$1e-6$ &$1e-3$ 	    & 64 &  5 &  74	 & 4.60e-3   & 8.46e-4   & 5.50e-4\\
	$1e-7$ &$3.16e-4$ 	& 88 &  8 &  104  & 1.04e-3   & 2.65e-4   & 2.14e-4\\
\bottomrule 
\end{tabular}
\caption[Space-time errors, coarse FEM]{Number of basis functions and space-time errors for the velocity and pressure. FEM solutions obtained on a coarse mesh, bypass application. Second column: number of selected velocity modes ($\# \uu$  basis). Third column: number of selected pressure modes ($\# p$  basis). Forth column: total number of basis functions composing the reduced space ($\#$ basis = $\# \uu$  basis + $2 \times \# p$  basis).}
\label{tab:6.1}
\end{table}

\begin{figure}[h]
\centering
% \subfloat[$tol =$ 1e-2]{
%	\includegraphics[width=5.0cm]{images/errRes1e_2}
% }
% \subfloat[$tol =$ 1e-3]{
%	\includegraphics[width=5.0cm]{images/errRes1e_3}
% }
% \subfloat[$tol =$ 1e-4]{
% 	\includegraphics[width=5.0cm]{images/errRes1e_4}
% }\\
 \subfloat[$tol =$ 1e-5]{
	\includegraphics[width=5.0cm]{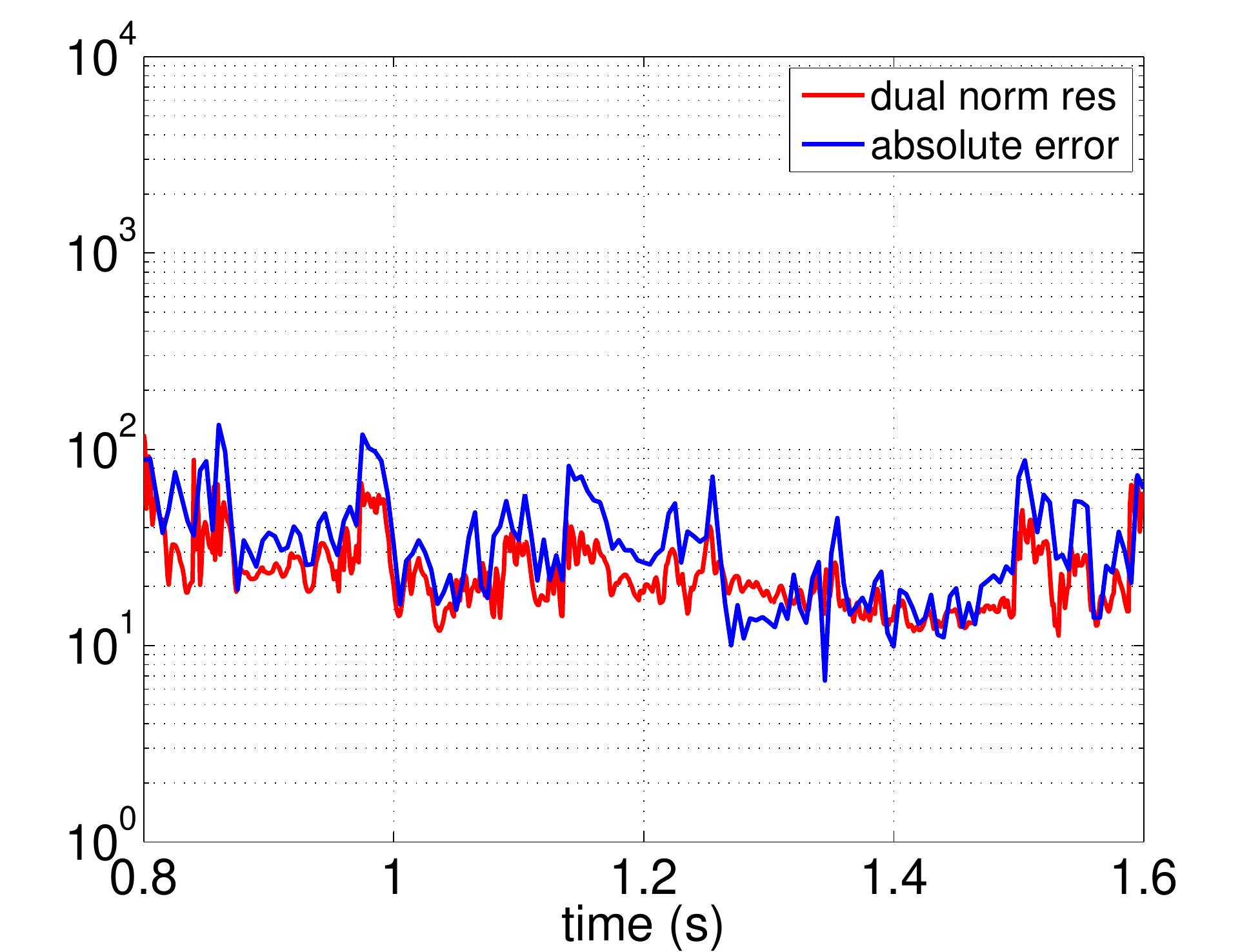}
 }
 \subfloat[$tol =$ 1e-6]{
	\includegraphics[width=5.0cm]{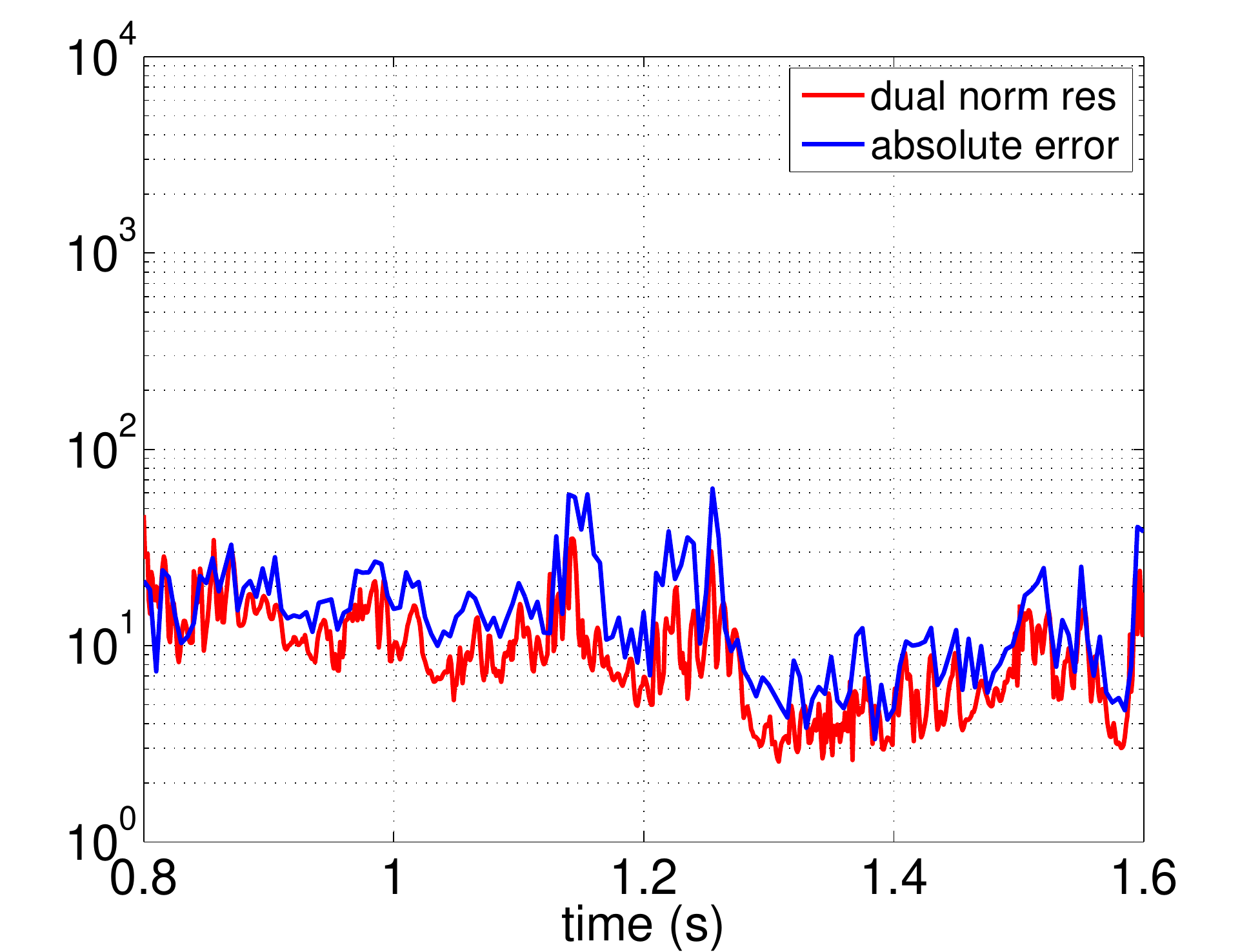}
 }
 \subfloat[$tol =$ 1e-7]{
	\includegraphics[width=5.0cm]{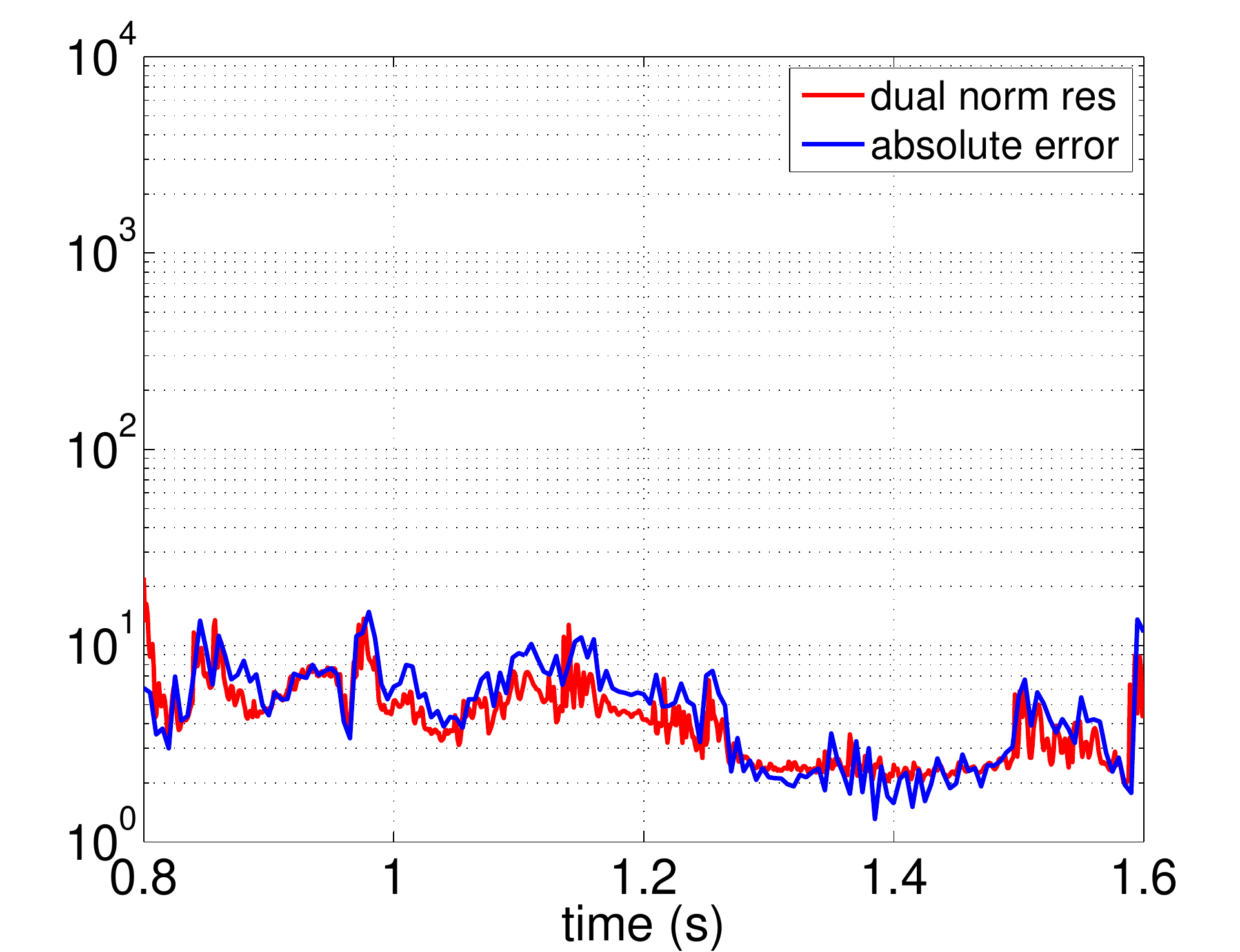}
}
\caption{Dual norms of the residuals and norms of the global errors with respect to time for different choices of the POD tolerance $tol$.}
\label{fig:6.2}
\end{figure}

\subsubsection{Application of the greedy enriched POD algorithm}
The POD algorithm provided satisfactory results and we were able to reduce the approximation space dimension from $10^5$ to $10-100$. In this section we aim at comparing the greedy enrichment algorithm with the POD one, in order to understand if using different basis functions than POD modes provides the same quality of reduced approximations. {\color{black}Thus, we compare the magnitude of the reduced approximation errors obtained using a reduced space generated through a standard POD algorithm with the ones obtained using the POD coupled with the greedy enrichment as introduced in Section \ref{sec:Greedy}.} We recall that the snapshots sample represents a subset of the time instants we solve in the unsteady simulation: indeed we store only the $25\%$ of the time instants solutions computed.
As there is no error bound available, we use the dual norm of the residual as surrogate. This is a rough approximation, also because the real error includes time integration, while the dual norm of the residual can only represent a space error. Of course we do not expect the greedy enriched POD to perform better, on the contrary it can have (and actually has) limitations. 

\textbf{Remark} {\color{black} We are interested to simulate a fluid-dynamics phenomena with cyclic inputs. Typically in haemodynamics applications, we are interested in
several heartbeats. Thus, instead of performing the greedy research only on one single heartbeat, we exploit as much as we can the \emph{information} on the truth solutions coming from the snapshots. For each single snapshot $\UUo^{n_s}_h$, we perform a simulation that starts from the initial time $t_{n_s}$ and ends at $t_{n_s + N_T} = t_{n_s} + 0.8$. We define a vector index $\mathbf n = (n_T, n_S)$ with $n_T=n_S+n$ such that $\UUo^{\mathbf n}_h = \UUo^{n_S, n_T}_h$ being the approximate solution at time $t^{n_T}$ obtained starting from the initial condition $\UUo^{n_s}_h$. We define the set of indexes $\mathcal N = \{(n_T,n_S) : n_T=n_S + n , n \in \mathcal N \text{ and } n_S \in \mathcal N_S\}$. The generalization of the greedy enrichment presented in Section \ref{sec:Greedy} is straightforwards substituting $n$ with $\mathbf n$. In particular, the selection of the \emph{worst} approximated index $n^*$ in the greedy enrichment can be generalized as follows:
\begin{equation*}
{\mathbf n}^* = \arg \max_{\mathbf n \in \mathcal N} \|r_N^{\mathbf n}(\mathbf W_h)\|_{\XX'}.
\end{equation*}}
%However, we study this way of enrichment because of its low costs
 
%\paragraph{Snapshots sample time step: 0.005}
In order to initialize the greedy enrichment algorithm we compute a POD basis fixing the tolerance $tol=1e-5$ (50 velocity modes, 3 pressure modes, 3 supremizers). To compare the POD approximation with the greedy enriched one, we augment the initial reduced space with two strategies. On one side, we apply the greedy enrichment and, at each iteration, we add the triplet of functions selected by the largest dual norm of the residual in space. On the other side, we augment the basis by adding, at each algorithm iteration, one POD mode for the velocity and one POD mode for the pressure with its associated supremizer. In both cases, at each iteration, we increase the reduced space dimensions of three units. The results obtained using only POD modes are displayed in black and addressed with the label \emph{POD}, while the results obtained with the greedy enrichment are shown in red and addressed with the label \emph{Greedy enriched POD} (see Figures \ref{fig:6.6}).
%\begin{figure}[t]
%\centering
%	\includegraphics[width=10cm]{images/snapSelection}
%\caption{Example of the time instants selection for which the high fidelity solution is stored.}
%\label{fig:5.0}
%\end{figure}

From Figures \ref{fig:6.6}, we note that the decrements of the errors in the greedy enrichment algorithm are slower than when adding POD modes. Nevertheless, we notice that both the space-time pressure error and residuals are comparable when adding POD modes or greedy basis functions (see Figure \ref{fig:6.6} (b) and (c)). On the contrary the decrements of the velocity is much slower when we use the greedy enrichment with respect to adding POD modes. We recall, however, than the residual is mostly related to the pressure error component.
%, in particular, we remark that the dual norm of the space-time residual in one case (see Figure \ref{fig:6.7} (a)) decreases faster using the greedy basis than using the POD modes. 
%Thus, we may assume that the greedy enrichment can be used when POD modes are not \emph{available}. This is the case, for example, when a coarse sample time step is chosen for the snapshot selection. We are going to detail this case in the next paragraph. Another possible application where it is interesting to test the greedy enrichment is represented by the case in which a parameter is introduced in the RFSI problem after POD generations and no high-fidelity solutions are available for different parameter realizations. In this situation using the greedy enrichment is much cheaper than generating new POD modes, that would require the computation of the finite element snapshots for new parameter realizations. We are going to detail this case in Section 6.2.

\begin{figure}[h]
\centering
 \subfloat[Space-time velocity error]{
	\includegraphics[width=5cm,height=4cm]{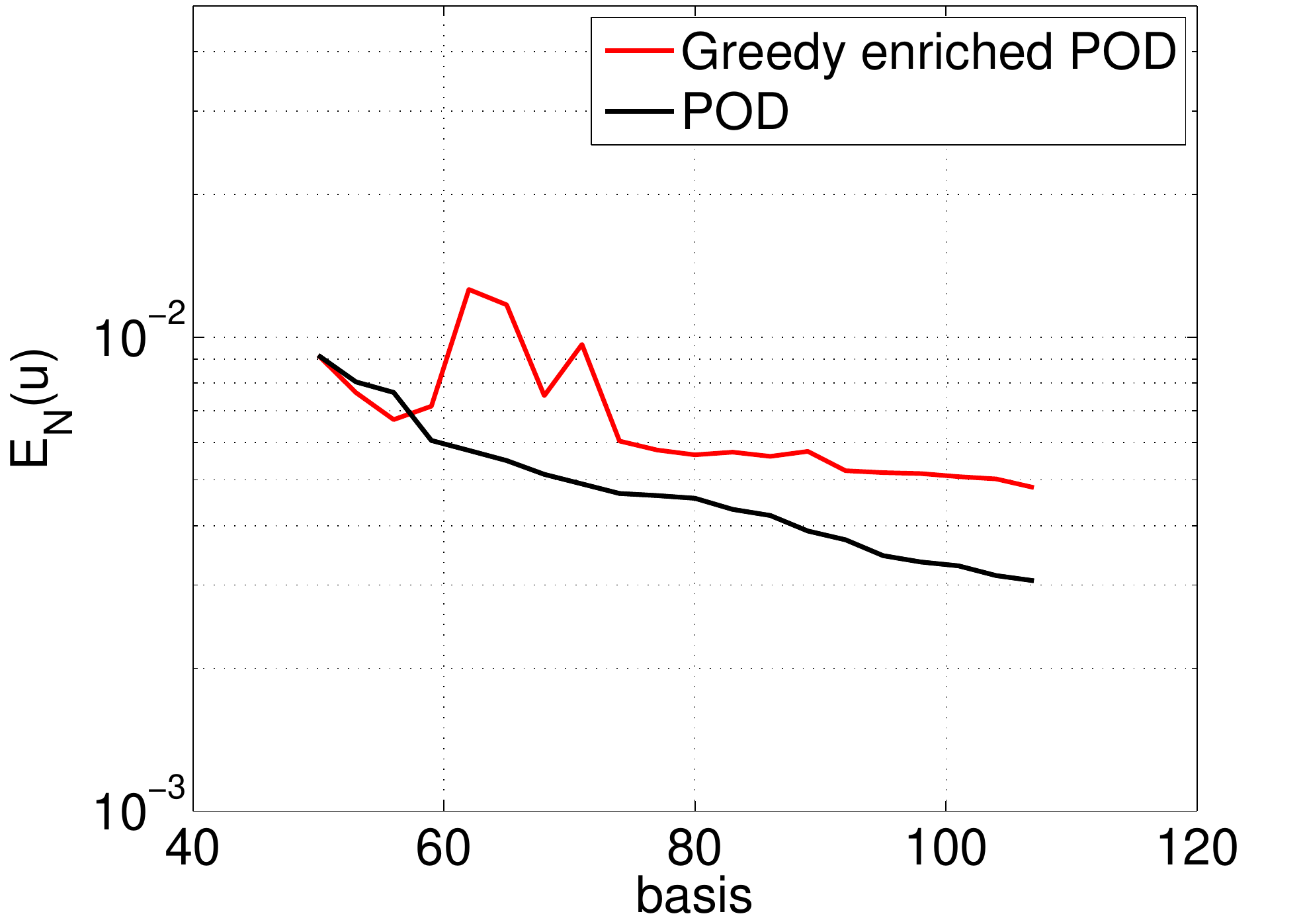}
 }
 \subfloat[Space-time pressure error]{
	\includegraphics[width=5cm,height=4.2cm]{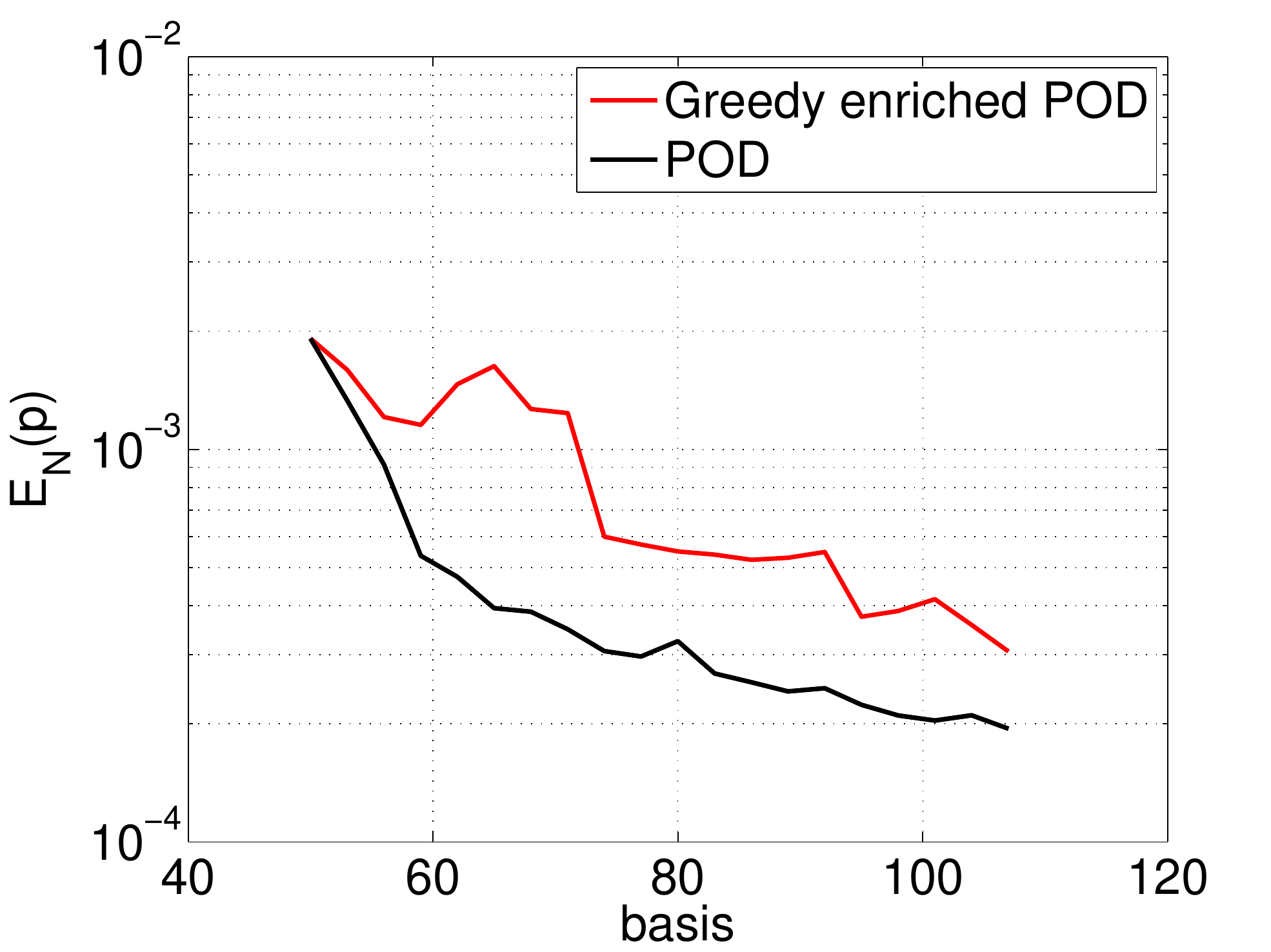}
 }
 \subfloat[Space-time dual norm of the residuals]{
	\includegraphics[width=5cm,height=4cm]{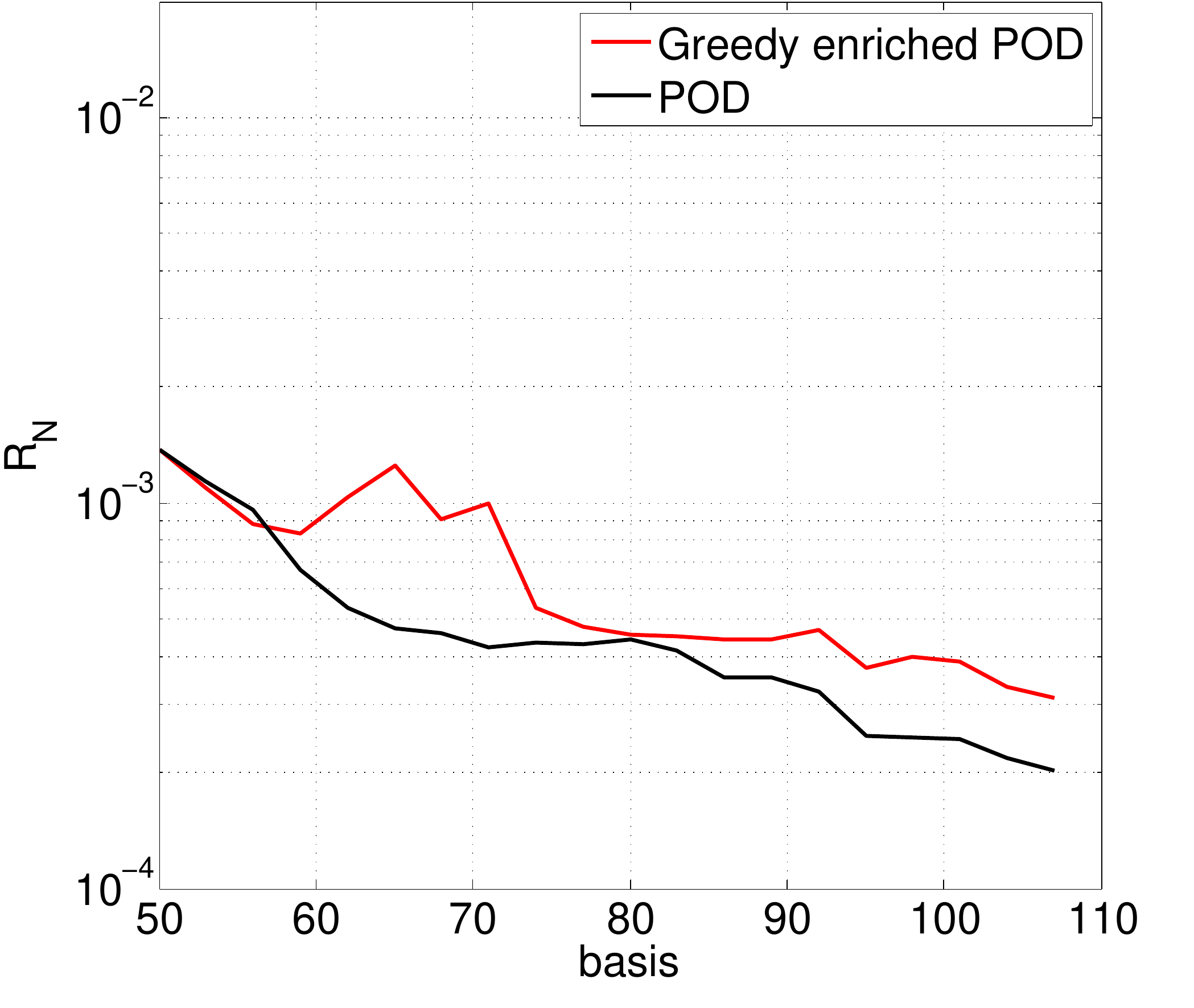}
	}
\caption{Space-time errors versus the number of basis, comparison between the standard POD and the greedy enriched POD, starting from 50 POD basis functions. FEM solutions obtained on a coarse mesh, bypass application.}
\label{fig:6.6}
\end{figure}

\subsection{Realistic case}
\subsubsection{Application of the POD algorithm}
In this section we perform a discretization reduction of the RFSI model applied to the femoropopliteal bypass case, where the high fidelity approximations are computed using a fine mesh. 
%The physical parameters, boundary conditions and the finite element discretization are the same as in Section 3.2.3. 
As before, a parabolic velocity profile is imposed at the inlet section and a mean pressure condition at the outlet. The $\mathbb P^1$+Bubble-$\mathbb P^1$ discretization yields 1'410'475 degrees of freedom on the fine mesh.
We first test the discretization reduction using a standard POD procedure: we compute the high fidelity numerical solutions for two heartbeats with a time step $\Delta t=0.001$ and we store the ones related to the second heartbeat every five time steps. Thus, $\mathcal N_T = {1,2,3,\ldots,800}$ and  $\mathcal N_S = {5,10,15,..,800}$.
%The snapshots are then represented by the following finite element solutions:
%\begin{equation*}
%\UU_h^{n_S} \approx \UU_h(t_{n_S}) \quad \text{ where } t_{n_0} = 0.805, t_{n_1}=0.810, t_{n_2}=0.815,t_{n_3}=0.820, .. , t_{N_S} = 1.6 \text{s}.
%\end{equation*}
We compute the Gramian matrices associated to the 160 snapshots $\UU_h^{n_k}$, separating the velocity and pressure components. We denote $\lambda^\uu_k$ and $\lambda^p_k$ for $k=1,..N_S$ the eigenvalues associated to the decomposition of the correlation matrices of the velocity and pressure, respectively (see Section \ref{sec:POD}). In both cases they decrease exponentially fast. The eigenvalues $\lambda^p_k$ associated to the pressure snapshots (Figure \ref{fig:6.2.0} (b)) decrease faster than the ones associated to the velocity (Figure \ref{fig:6.2.0}(a)). Thus, by fixing the same tolerance, we expect that a fewer number of pressure modes will be selected with respect to the velocity ones.

\begin{figure}[t]
\centering
 \subfloat[Velocity eigenvalues, $\lambda^\uu_k$]{
	\includegraphics[width=5.4cm]{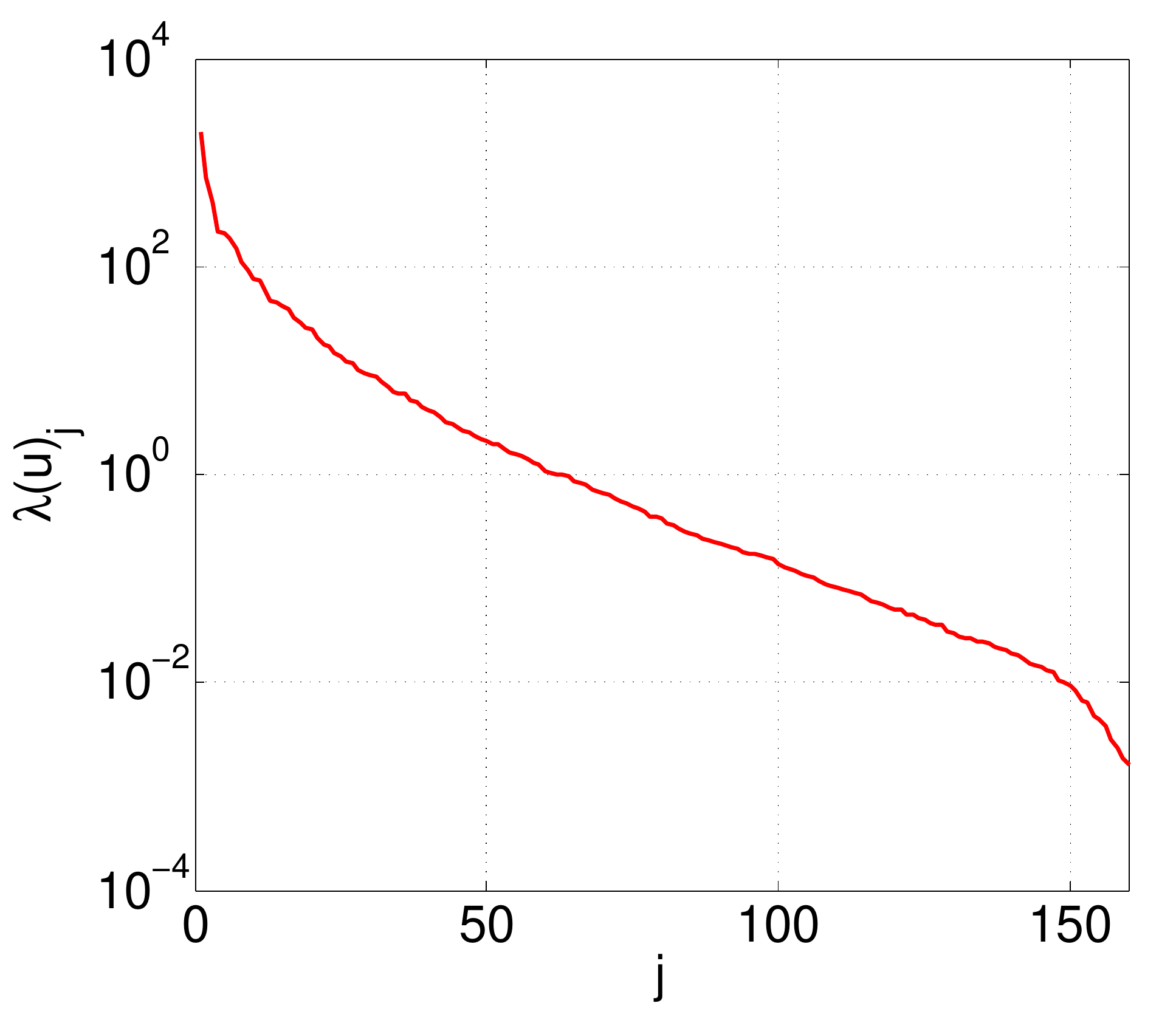}
 }
 \subfloat[Pressure eigenvalues, $\lambda^p_k$]{
	\includegraphics[width=5.6cm]{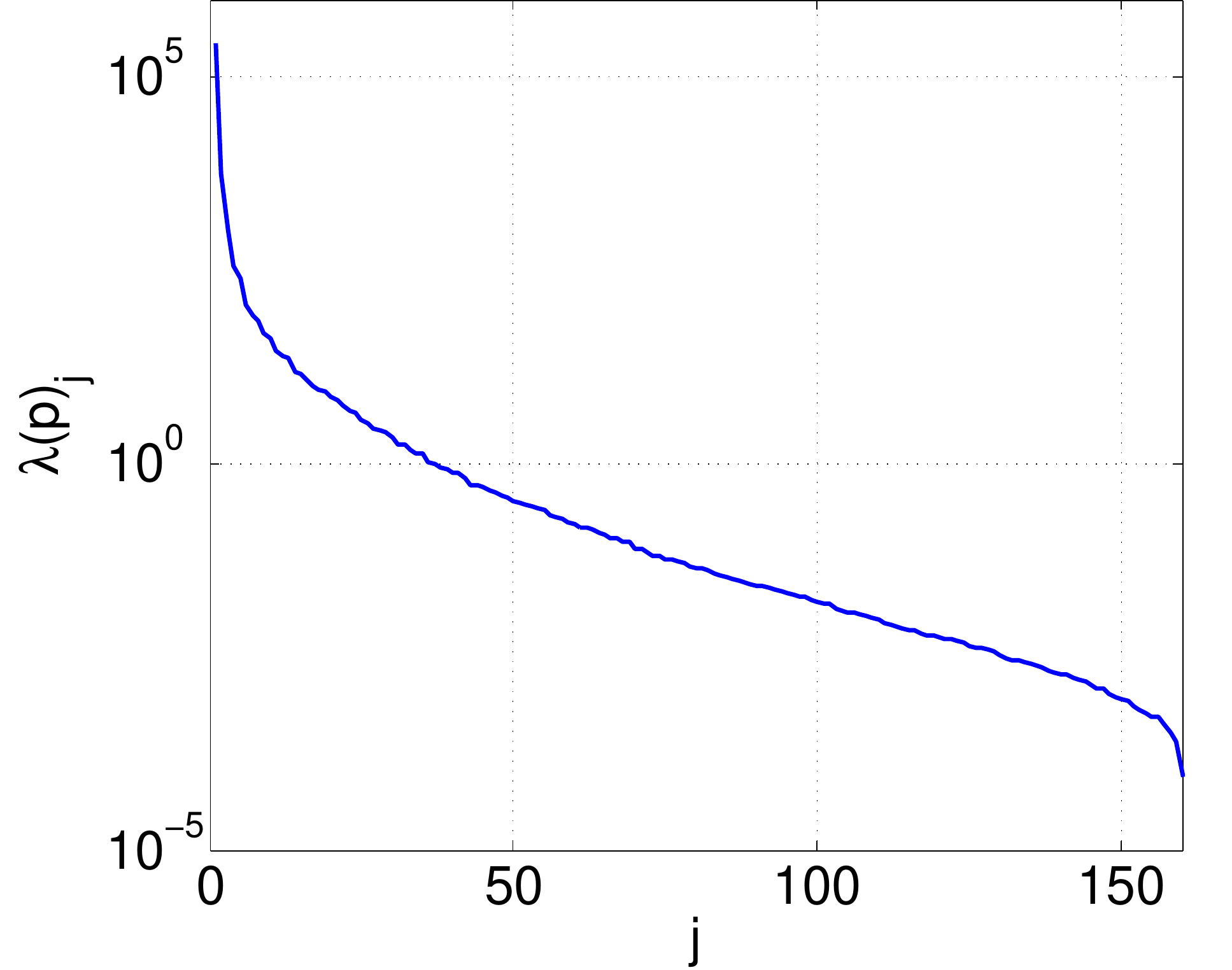}
 }
\caption{Velocity and pressure eigenvalues computed with 160 snapshots sampled every 0.005 s.}
\label{fig:6.2.0}
\end{figure}

We compute the POD reduced spaces using three different values for the tolerance: $tol \in \{1e-4, 1e-5, 1e-6\}$. As it was done in Section \ref{sec:TestCase}, in Table \ref{tab:6.2.1}, we record the selected number of modes and we compute the space-time errors $E_N(\uu)$ and $E_N(p)$ of the velocity and pressure, respectively. By taking advantage of the generated reduced space, at each time iteration we solve the reduced system and we compute a linear functional of the approximate solution that evaluates the outlet flow rate. We record the computational time associated with the offline and online computations in Table \ref{tab:6.2.0}. We can appreciate that the resolution of the reduced problem combined with the evaluation of a linear output functional is performed in almost real time: using 35 basis functions we solve 1.6 \emph{physical} seconds in 0.8 \emph{computational} seconds, while a ten heartbeats simulation (8 physical seconds) selecting the POD space with 54 basis functions takes 12.6 s on a notebook. Performing the same simulation with the high fidelity model would have taken around 40 hours on 256 processors of a supercomputer.
The offline costs of the POD reduction (without the snapshots generation) are reported in the third column of Table \ref{tab:6.2.0}. We remark that most of this time is spent in the generation of the structures for the residual evaluation.

Note that the POD model errors $E_N(\uu)$ and $E_N(p)$ decrease significantly when increasing the number of basis functions, as it is reported from both the values of Table \ref{tab:6.2.1}. Once again, we notice that the dual norm of the residual $\|r_N^{n+1}(\mathbf W_h)\|_{\XX'}$ is a good indicator of the approximation error $\|\UU_N^{n_k} - \UU_h^{n_k}\|_{\XX}$ (see Figure \ref{fig:6.2.3}).

{\color{black} In the femoropopliteal bypass application, we are interested in measuring also the errors on the output of interests. Being $\sigma^{n_S}$ the stress tensor and $\bm n$ the normal vector to the surface $\Gamma$, we compute the wall shear stress as $\bm \tau^{n_S}:=\sigma^{n_S}\bm n - (\sigma^{n_S}\bm n \cdot \bm n)\bm n$ and we also consider the averaged wall shear stresses on a generic area $A$:  $\tau_A^{n_S} = 1/A \int_A |\bm \tau^{n_S}|dA$.
%\begin{equation*}
%wss_i^{n_k}: = \dfrac{1}{A_i}\int_{A_i} |\mathbf t^{n_k} - (\mathbf t^{n_k} \cdot \mathbf n) \mathbf n| dA_i \qquad \text{ for } i = 1,2,3, 
%\end{equation*}
%being $A=A_1$ or $A=A_2$ and $A_3$ the areas defined in Figure \ref{fig:6.2.15}(a). 
%%We measure the correspondent %wall shear stress errors: the relative error of the averaged wall shear stress at time $t_{n_k}$:
%%\begin{equation}
%%\varepsilon_N(wss_i^{n_k}):=\dfrac{|wss_{i,N}^{n_k} - wss_{i,h}^{n_k}|}{wss_{i,h}};
%%\end{equation}
%space-time error of the averaged wall shear stress
%\begin{equation*}
%E_N(\mathbf tau_{A}):=\dfrac{\bigg(\sum_{n_S \in \mathcal N_S} |\mathbf tau_{A,N}^{n_S} - \mathbf tau_{A,h}^{n_S}|^2\bigg)^{1/2}}{\bigg(\sum_{n_S \in \mathcal N_S}\big(\mathbf tau_{A,h}^{n_S}\big)^2\bigg)^{1/2}} \qquad \text{ for } A = A_1,A_2.
%\end{equation*}
%We note that the accuracy in the approximation of the wall shear stress increases when considering a higher number of basis (see Table \ref{fig:6.2.3}). Moreover, the space-time errors $E_N(\mathbf \tau_A)$ (see Table \ref{fig:6.2.3}) are of the same order of magnitude of the velocity errors (see Table \ref{fig:6.2.1}). 
We remark that to properly estimate the selected output of interest we need accurate high fidelity solutions with a mesh refined at the wall \cite{Marchandise:2011aa}. Reducing the dimension of the finite element space does not lead to the same results that we obtain reducing the degrees of freedom using the POD decomposition. In fact, the wall shear stress values computed using a coarse finite element space underestimate considerably the values obtained with the fine grid (see Figure \ref{fig:6.2.15}), while the results obtained with the POD reduced approximation mostly overlapped with the ones computed with the finite element discretization.
}
\begin{table}[H]
\centering
\begin{tabular}{ccccccccc}
\toprule
	$tol$ & $\# \uu$  basis 	& $\# p$ basis & $\#$ tot basis & $E_N(\uu)$ & $E_N(p)$ & $R_N^{\UU}$  \\ 
\midrule
	$1e-4$ 	& 31 & 2 & 35 & 5.505e-2  & 1.188e-2  & 8.599e-3\\
	$1e-5$ 	& 48 & 3 & 54 & 9.840e-3  & 1.910e-3  & 1.441e-3 \\
	$1e-6$ 	& 68 & 5 & 78 & 5.074e-3  & 9.131e-4  & 6.264e-4 \\
\bottomrule 
\end{tabular}
\caption[Space-time errors, fine FEM]{Number of basis functions and space-time errors for the velocity and pressure. FEM solutions obtained on a fine mesh, bypass application.}
\label{tab:6.2.1}
\end{table}

\begin{table}[H]
\centering
\begin{tabular}{ccccc}
\toprule
	$tol$ & $\#$ tot basis & CPU time $X_N^{POD}$ & CPU time 2HB - RB &  CPU time 2HB - FE  \\ 
\midrule
	$1e-4$  & 35 & $\sim$ 38 min & 0.84 s & $\approx$  28800 s (8 hrs)\\
	$1e-5$ & 54  & $\sim$ 85 min & 2.49 s & $\approx$ 28800 s (8 hrs)\\
	$1e-6$ & 78  & $\sim$ 172 min & 6.84 s& $\approx$ 28800 s (8 hrs)\\
\bottomrule 
\end{tabular}
\caption[Computational times]{CPU time $X_N^{POD}$: offline computations costs for the generation of the POD reduced spaces (without the finite element computations) on 512 processors; CPU time 2HB - RB: online computational time corresponding to the simulation of 2 heartbeats (2HB) on a personal laptop; CPU time 2HB - FE: finite element computational time corresponding to the simulation of 2 heartbeats (2HB) on 256 processes on a supercomputer. }
\label{tab:6.2.0}
\end{table}

\begin{figure}[H]
\centering
 \subfloat[$tol =$ 1e-4]{
 	\includegraphics[width=5.0cm]{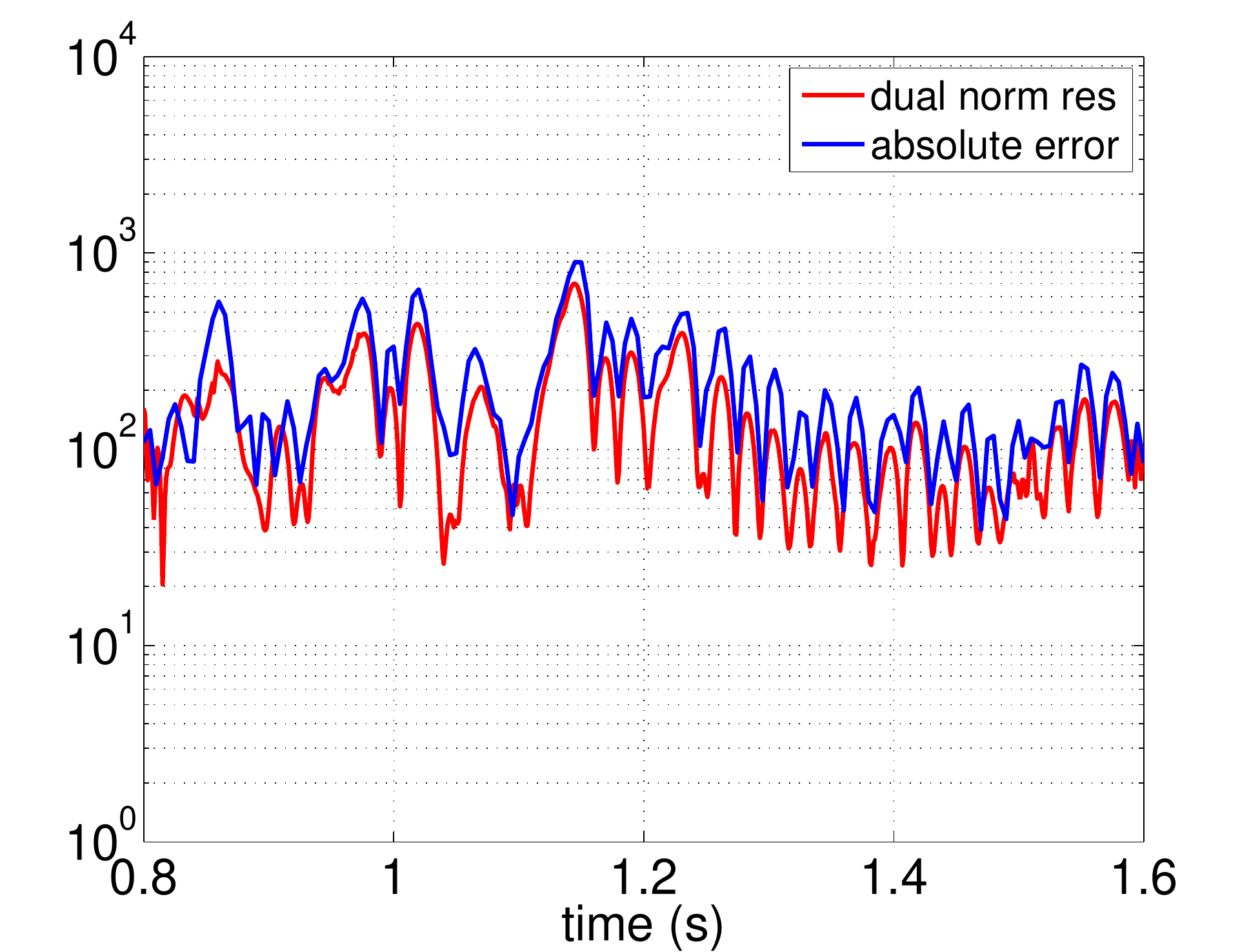}
 }
 \subfloat[$tol =$ 1e-5]{
	\includegraphics[width=5.0cm]{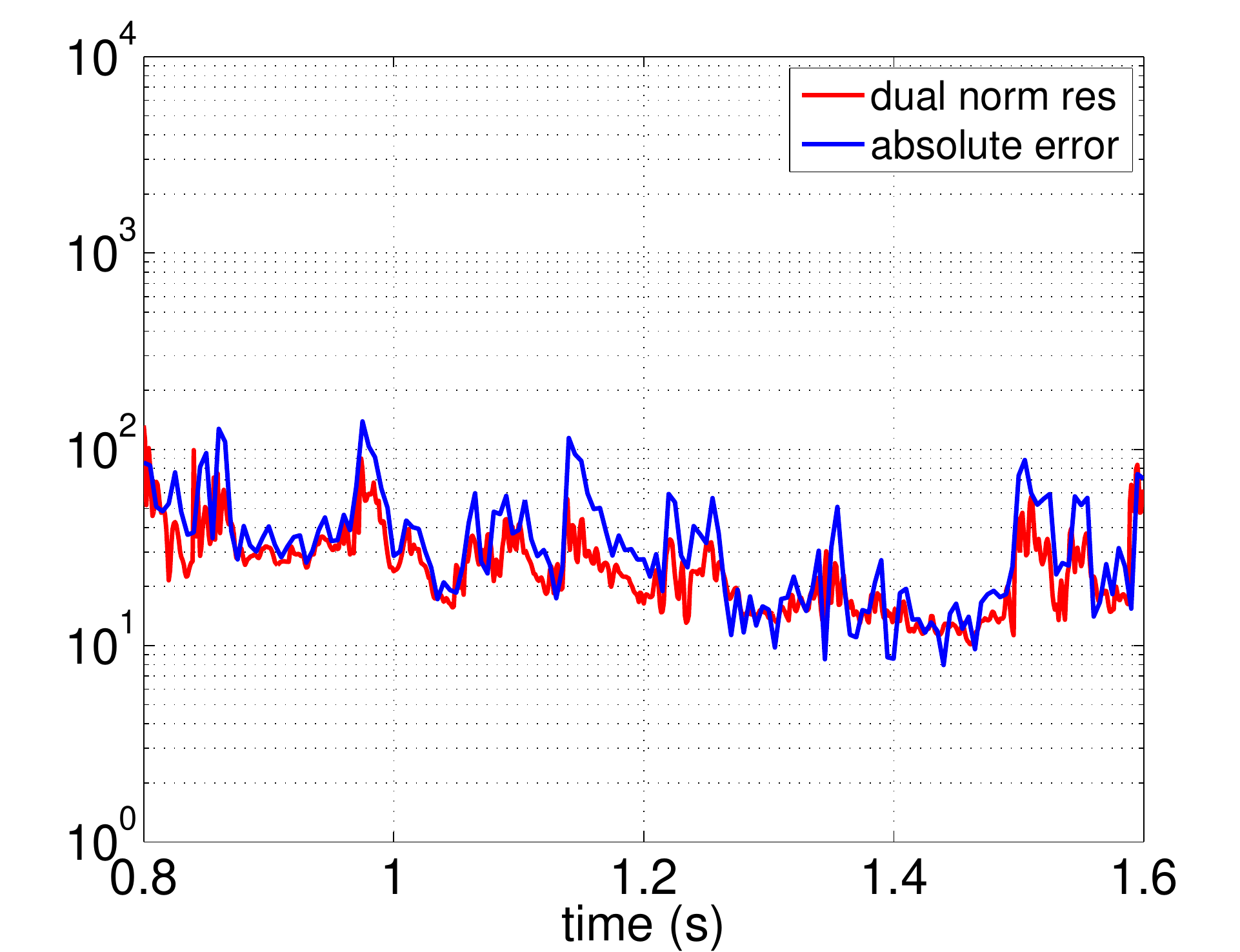}
 }
 \subfloat[$tol =$ 1e-6]{
	\includegraphics[width=5.0cm]{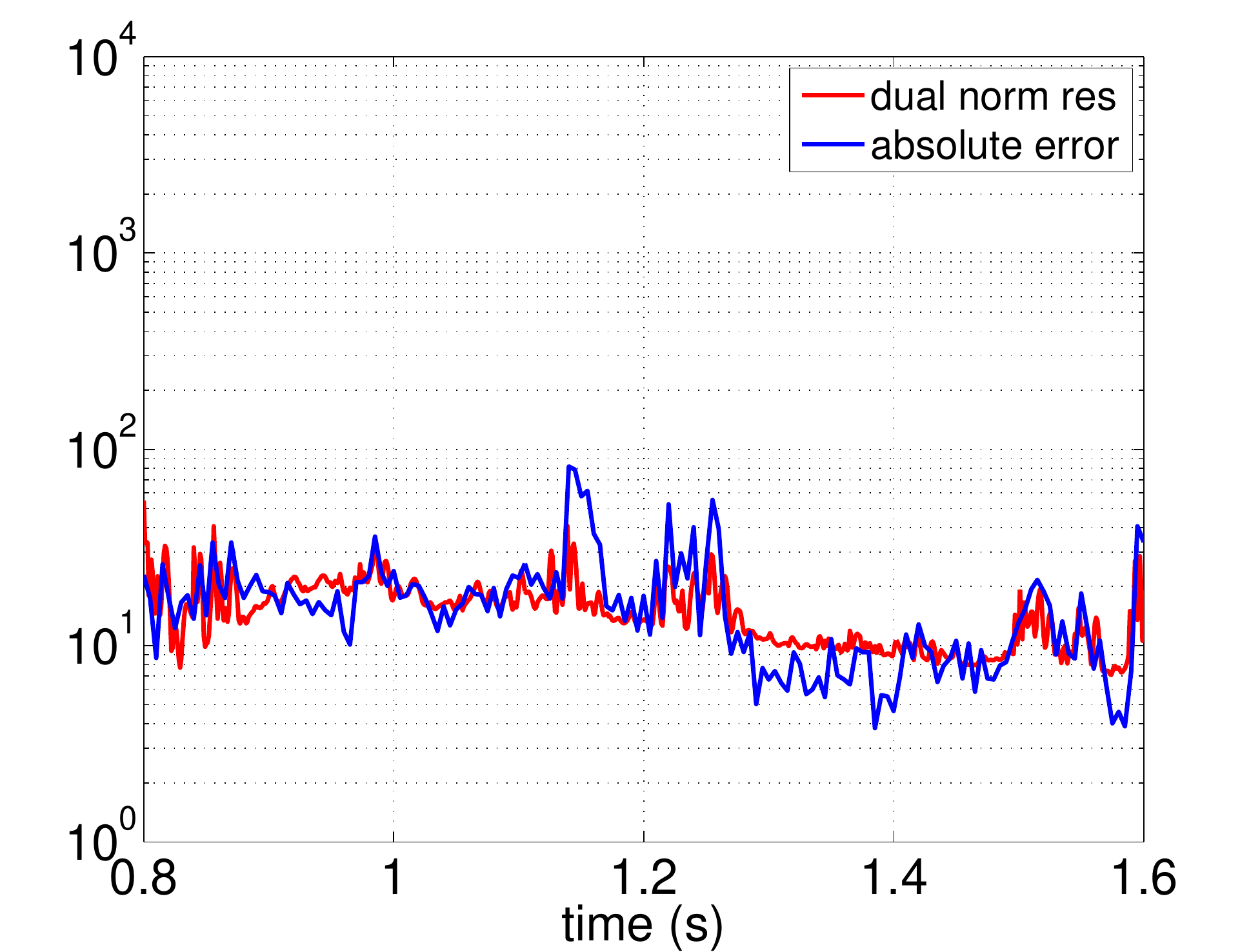}
 }
\caption{Dual norm of the residuals and $X$-norm of the errors with respect to the time instant for different choices of the POD tolerance $tol$. FEM solutions obtained on a fine mesh, bypass application.}
\label{fig:6.2.3}
\end{figure}

\begin{figure}[H]
\centering
  \subfloat[Areas location]{
	\includegraphics[width=2.0cm]
	{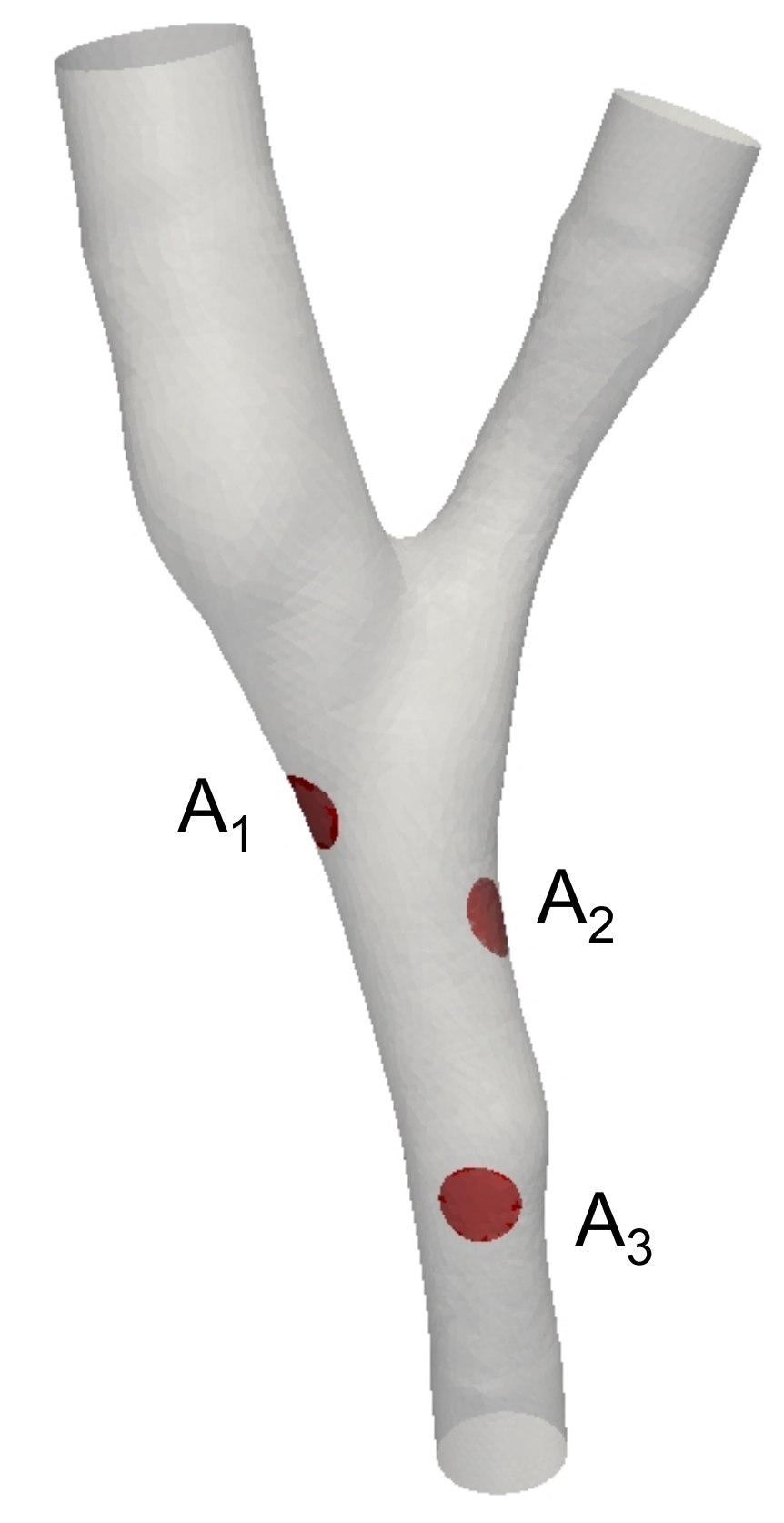}
 }
 \subfloat[Wss absolute values - Area 1]{
 	\includegraphics[width=5cm]{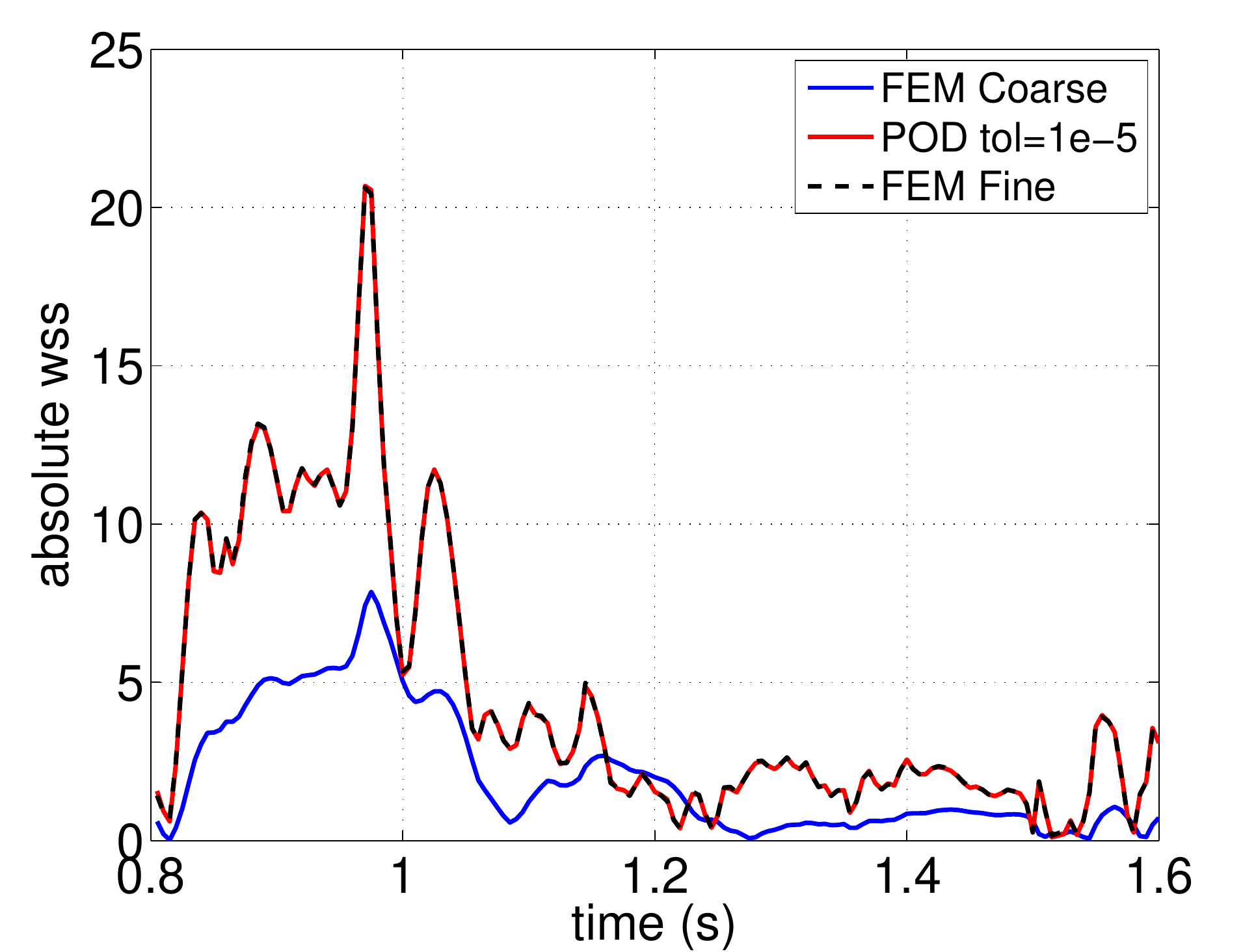}
 }
 \subfloat[Wss absolute values - Area 2]{
	\includegraphics[width=5cm]{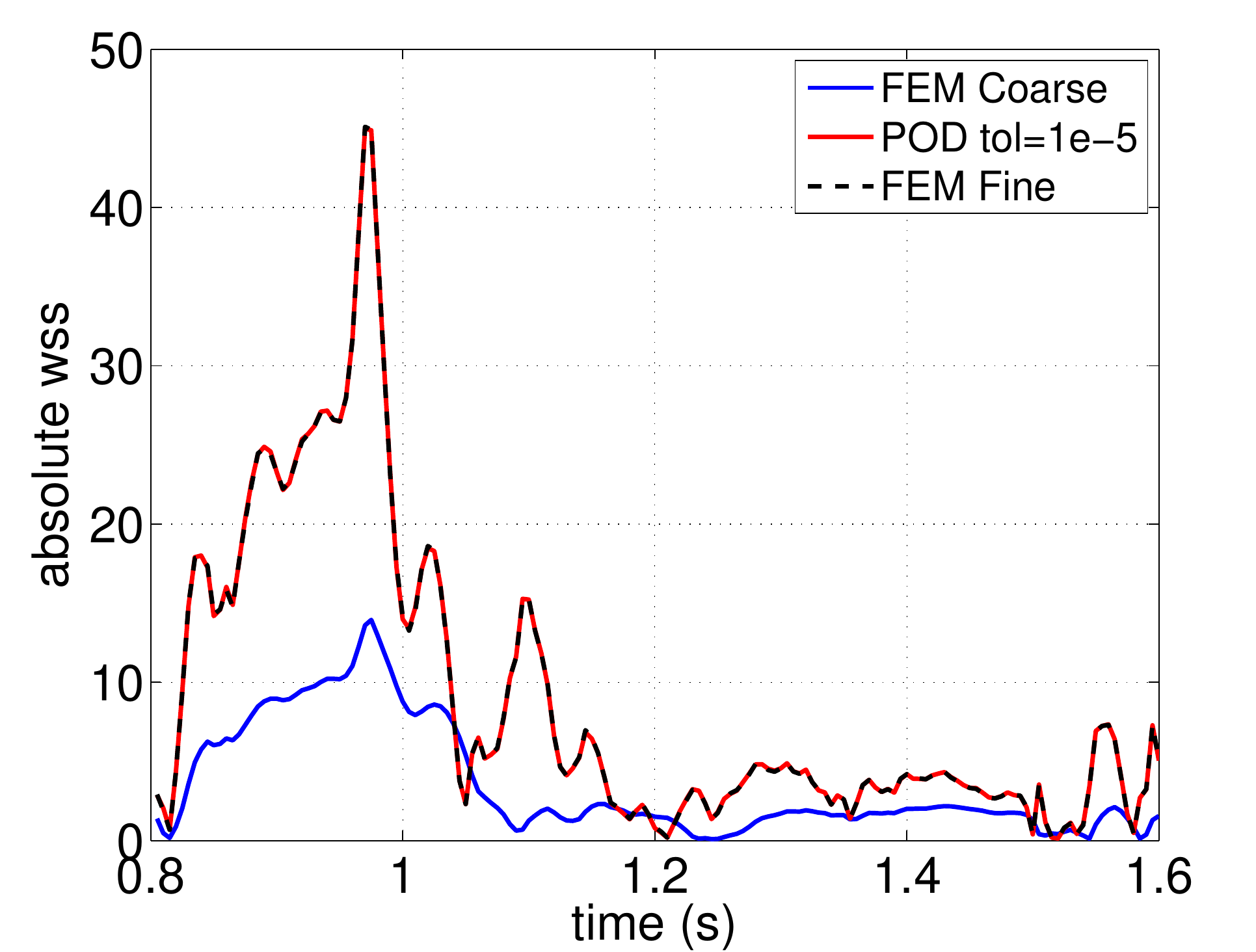}
}
\caption{Averaged wall shear stress computed in different location of the interface. Comparisons between the values obtained with the finite element discretization on the fine mesh (FEM Fine), the coarse one (FEM Coarse) and the reduced basis approximation (POD tol = 1e-5)}
\label{fig:6.2.15}
\end{figure}

\subsubsection{Application of the greedy enriched algorithm with perturbed data}
In this section we apply the greedy enrichment in the case of perturbed boundary data. In particular, as in \eqref{eq:PerturbG}, we introduce a parameter in the inlet flow rate function representing a small perturbation with respect to a reference value. The \emph{perturbation} function $\theta(\alpha, t)$ is define as follows:
\begin{equation*}
\theta(\alpha,t) = 1 + \alpha \sin\bigg(\frac{2 \pi t}{0.8}\bigg)
\end{equation*}
where $\alpha$ is supposed to vary between 0 and 0.2. Thus, the maximum relative difference with the original flow rate is equal to the 20$ \%$. We denote with $\UUo^{\mathbf n}_*(\alpha)$ with $*=\{h\},\{N\}$ or $\{N,h\}$ the numerical solutions at the time instant $t^{\mathbf{n}}$ that depend on the parameter $\alpha$ and with $r_N^{\mathbf n}(\mathbf W_h; \alpha)$ the residuals.  
In the perturbed case, the algorithm steps in Section \ref{sec:Greedy} are modified as follows. First we perform a POD algorithm fixing $\alpha=\alpha_1=0.0$; the resulting reduced space is addressed with $\XX_N^{\alpha_1,POD}$, where the apex $\alpha_1$ denotes the choice of the $\alpha$ parameter. Then, we set $\alpha_2=0.2$:
\begin{itemize}
\item[1.] Generate the reduced basis solutions $\UUo_N^{\mathbf n}(\alpha_2)$, $\mathbf n \in \mathcal N$ by solving the reduced order problem \eqref{eq:5.36}.
\item[2.] Compute the dual norms of the residuals $\|r_N^{\mathbf n}(\mathbf W_h; \alpha_2)\|_{\XX'}$, $\mathbf n \in \mathcal N$, which are used as error indicators. 
\item[3.] Select $\mathbf{n}^*$ such that $\mathbf n^* = \arg \max_{\mathbf n \in \mathcal N} \|r_N^{\mathbf n}(\mathbf W_h; \alpha_2)\|_{\XX'}$.
\item[4.] Compute the $\UU_{N,h}^{\mathbf n^*}(\alpha_2)$ by solving the reduced order problem  \eqref{eq:5.37}.
\item[5.] Split $\UU_{N,h}^{\mathbf n^*}(\alpha_2)$ into its velocity and pressure components, $\uu_h^{\mathbf n^*}$ and $p_h^{\mathbf n^*}$, respectively. Compute the supremizer $\bm \sigma_{\mathbf n^*}$ associated with the pressure component.
\item [6.] Compute $\bm \phi^{\uu}$ representing the orthonormalization of the velocity function $\uu_h^{\mathbf n^*}$ with respect to the reduced space $\XX_N$, obtained with a Gram-Schmidt algorithm considering the scalar product $(\cdot,\cdot)_{\VV}$; similarly for $\phi^{p}$ and $p_h^{\mathbf n^*}$.
\item [7.] Build $\XX_{N+2}=\XX_{N}\oplus \{\bm \psi^{\uu},\bm \psi^{p}\}$ built as is \eqref{eq:finalBasis}.
\item [8.] Compute $\bm \phi^{\bm \sigma}$ representing the orthonormalization of the velocity function $\bm \sigma^{n^*}$ with respect to the reduced space $\XX_{N+2}$, obtained with a Gram-Schmidt algorithm considering the scalar product $(\cdot,\cdot)_{\VV}$.
%= \dfrac{\uu_h^{n^*} - \Pi_{X_N}\uu_h^{n^*}}{\|\uu_h^{n^*} - \Pi_{X_N}\uu_h^{n^*}\|_{\VV}}$ and similarly for $\phi^{p}$ and $\bm \phi^{\bm \sigma}$
\item [7.] Build $\XX_{N+3}=\XX_{N+2}\oplus \{\bm \psi^{\bm \sigma}\}$ built as is \eqref{eq:finalBasis}.
\item[8.] Update the structures for the online computation of the reduced solutions and the dual norms of the residuals.
\item[9.] Set $N = N+3$ and $\XX_N = \XX_{N+3}$. Repeat until a predefined stopping criterion is satisfied.
\end{itemize}
{\color{black} The real modification is indeed related to the fact that the initial POD is computed for $\alpha=\alpha_1=0$, while the greedy enrichment is performed fixing $\alpha=\alpha_2=0.2$. The resulting reduced space $\XX_N$ aims to represent a suitable space of approximation for both values of $\alpha$.} In the parametrized case, by using greedy enrichment we aim at saving a part of the offline computational costs: indeed, in a standard POD-Greedy procedure (see \cite{Haasdonk:2008aa}), each new evaluation of the parameter $\alpha$ requires the computation of the associated finite element solutions for each time instant $n \in \mathcal N_T$. In our application, this would require about 8 hours on 256 processors. Instead, during the greedy enrichment, we perform only one finite element resolution for a single time step, while the remaining computations are dedicated to reduced basis structures.

We test the greedy enrichment algorithm by initializing it with two different starting POD reduced spaces: in one case we consider the modes selected with $tol=1e-4$ (35 POD basis functions) and in the other one we consider the POD modes corresponding to $tol=1e-5$ (54 POD basis functions). In the first case, we enrich the space $\XX_{35}^{\alpha_1,POD}$ by adding 8 triplets selected by the greedy algorithm; we obtain the reduced space $\XX_{59}$. In the second case, starting from $\XX_{54}^{\alpha_1,POD}$, we enrich the space adding 12 triplets, obtaining $\XX_{90}$. All the errors and residuals computed and shown below are referred to the solutions obtained with $\alpha_2=0.2$. In particular, in Table \ref{tab:6.2.2} we report the velocity and pressure errors generated by the greedy enriched reduced spaces as well as the ones obtained with the standard POD ones. Moreover, we compute the space-time dual norm of the residual, scaled by the solution norm (sixth column of Table \ref{tab:6.2.2}). 

We note that the space-time velocity error does not decrease significantly neither when adding greedy basis functions nor when augmenting the number of selected POD modes. 
If we look at the pressure, using the greedy enrichment we manage to decrease its error more than if we use POD modes. Also the space-time dual norm of the residual is smaller when considering the greedy enriched space than the POD ones.

Regarding the offline costs, to generate the space $\XX_{59}$ starting form the $\XX_{35}^{0,POD}$, we perform 8 iterations of the greedy enrichment algorithm: this takes 82 minutes on 512 processors where the most of the time is devoted to the generation of the reduced structures for the residual evaluation. We remark that computing a standard POD reduced space for the parameter evaluation corresponding to $\alpha=0.2$ would require about 8 hours on 256 processors for the finite element computation of two periods, plus about 1 hour on 512 processors for the generation of the reduced space itself. 
\begin{table}[t]
\centering
\begin{tabular}{cccccc}
\toprule
   $\# \uu$  basis (Greedy) 	& $\# p$ basis (Greedy) & $\#$ tot basis & $E_N(\uu)$ & $E_N(p)$  & $R_N^{\UU}$  \\ 
\midrule
31 (0) & 2  (0) & 35  & 1.403e-01  &  1.256e-02 &  1.284e-02   \\
48 (0) & 3  (0) & 54  & 1.301e-01  &  3.878e-03 &  7.125e-03   \\
39 (8) & 10 (8) & 59 & 1.326e-01 &  2.966e-03   &  3.048e-03\\
68 (0) & 5  (0) & 78  & 1.217e-01  &  3.239e-03  &  6.594e-03   \\
60 (12) & 15 (12) & 90  & 1.261e-01  &  1.748e-03 &  1.630e-03   \\ 
\bottomrule 
\end{tabular}
\caption[Space-time residual, fine FEM]{Number of basis functions and space-time errors for the velocity and pressure. Femoropopliteal bypass application in which the high fidelity solutions are obtained using a finite element approximation on a fine mesh.}
\label{tab:6.2.2}
\end{table}

\begin{figure}[h!]
\centering
% \subfloat[35 Basis - POD]{
%	\includegraphics[width=5.0cm]{images/errResF2_1}
% }
 \subfloat[54 Basis - POD]{
	\includegraphics[width=5.0cm]{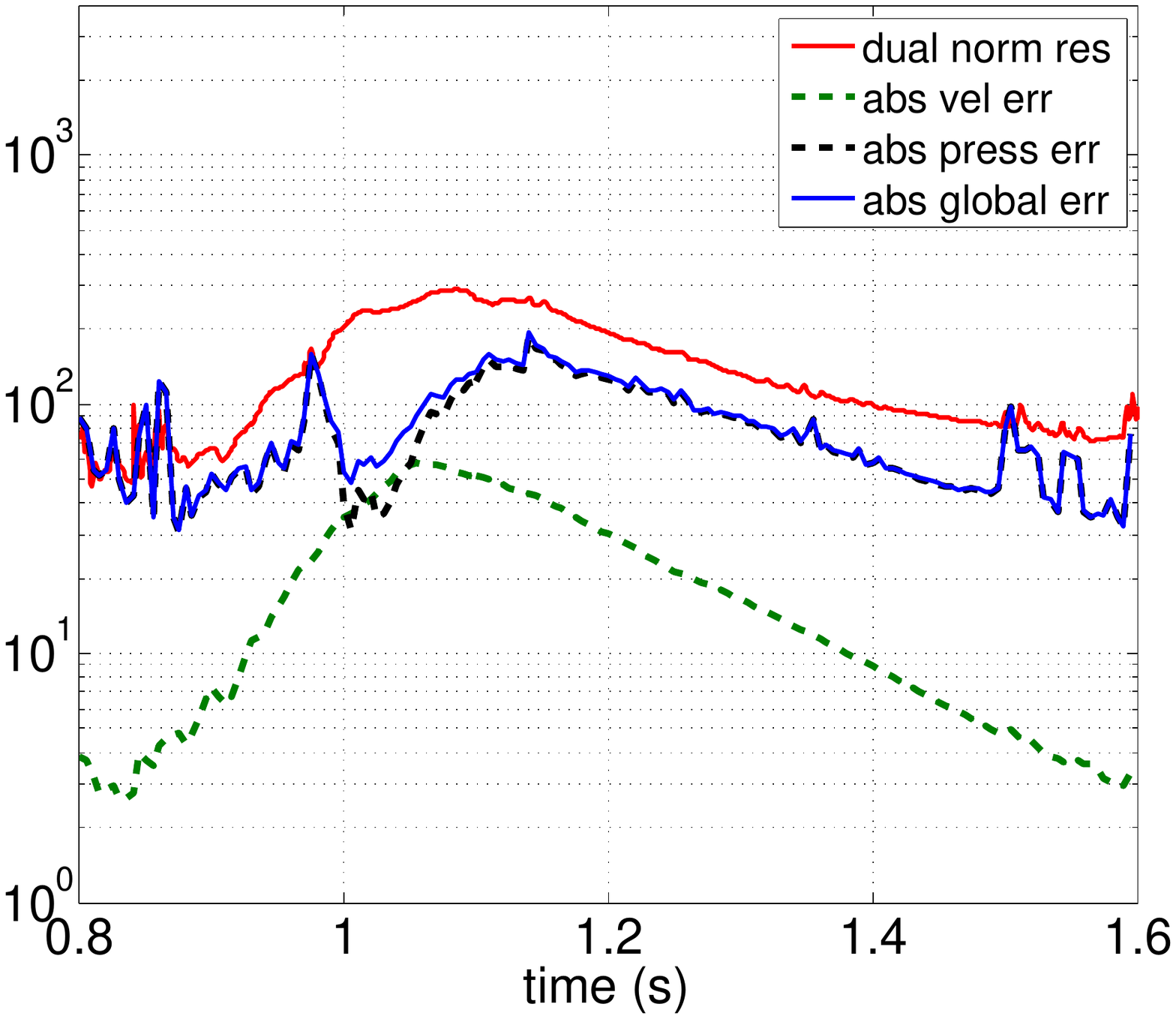}
 }
 \subfloat[78 Basis - POD]{
 	\includegraphics[width=5.0cm]{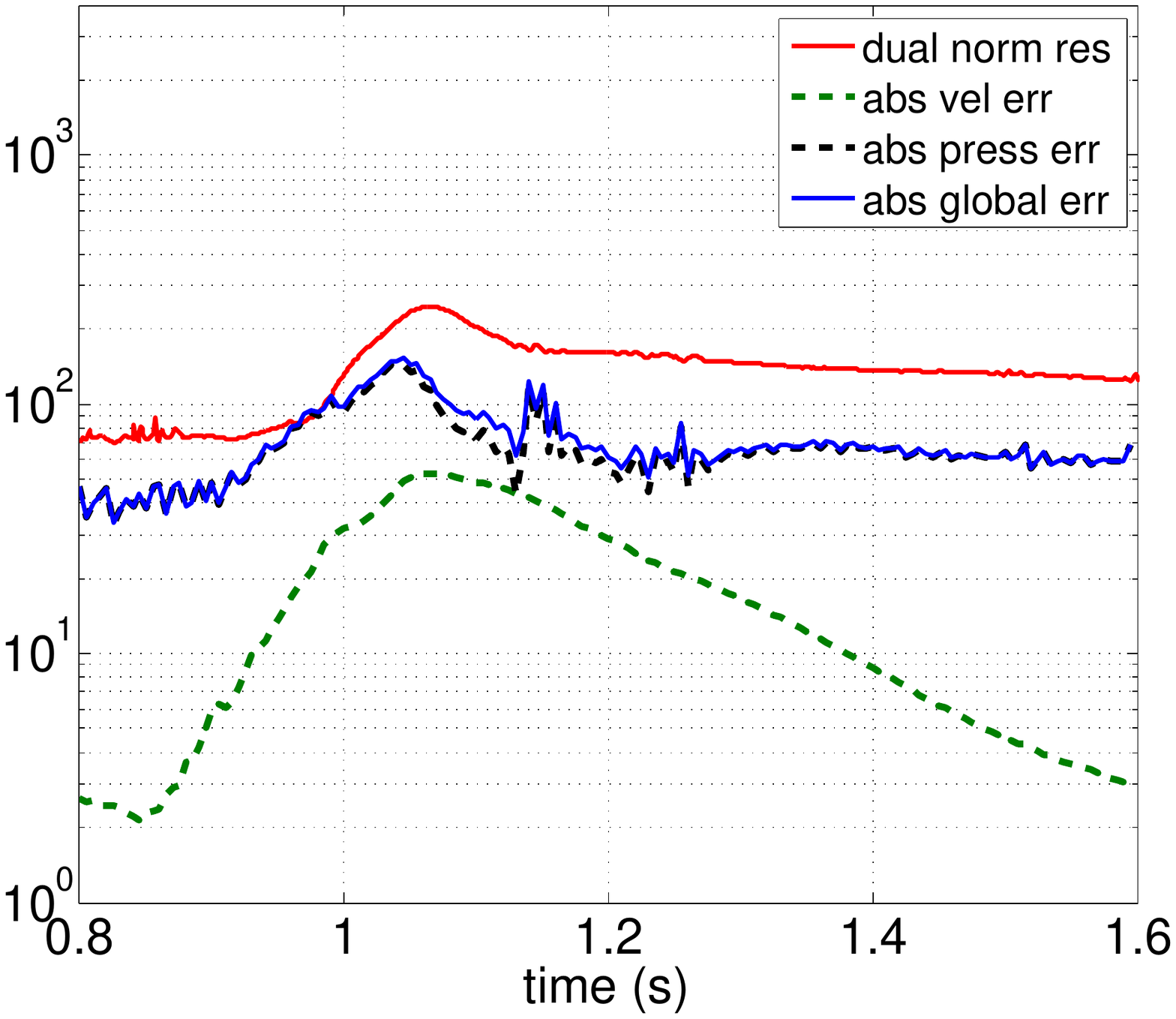}
 }
% \subfloat[59 Basis - Greedy enriched POD]{
%	\includegraphics[width=5.0cm]{images/errResF2_4}
% }
 \subfloat[90 Basis - Greedy enriched POD]{
	\includegraphics[width=5.0cm]{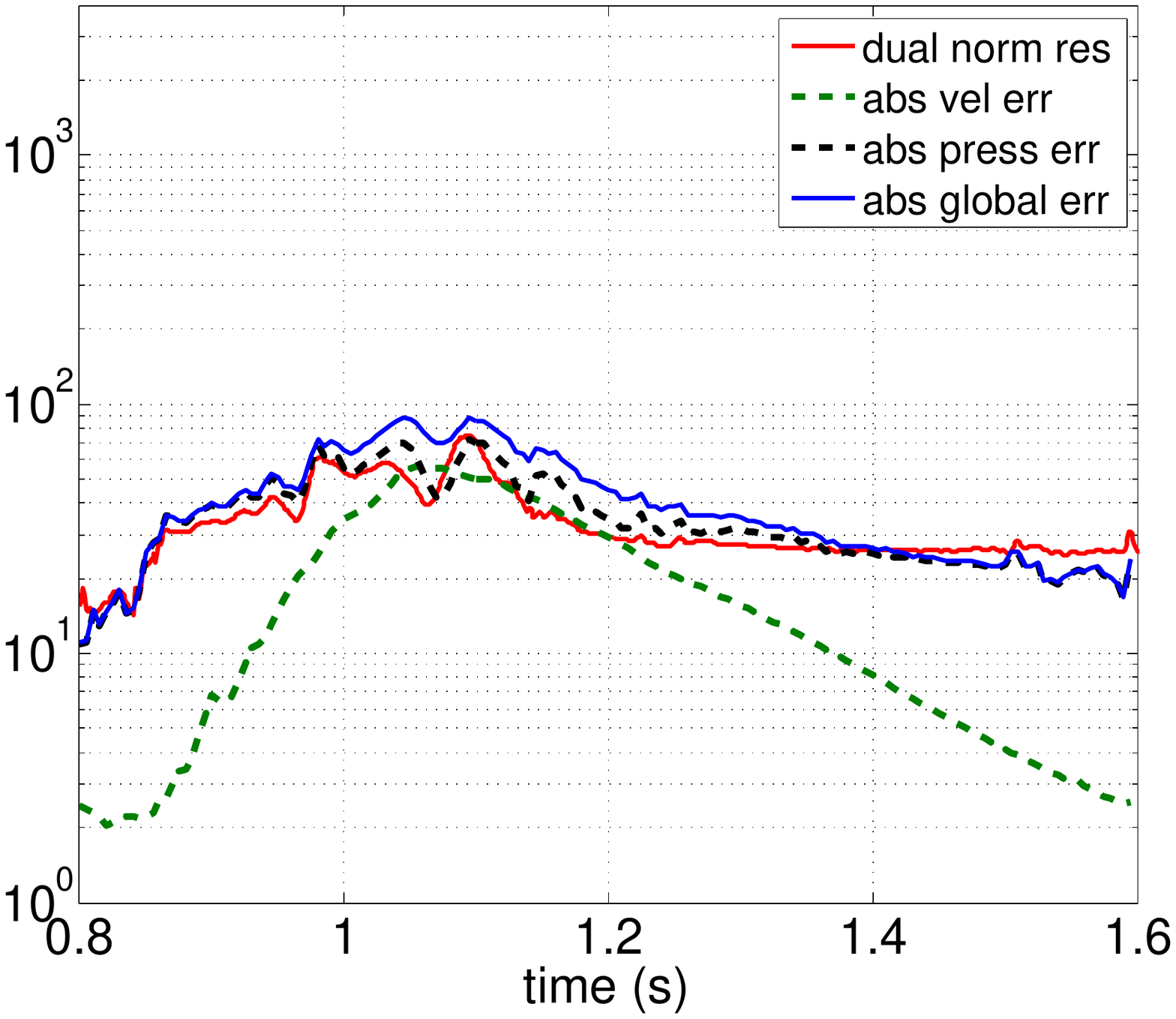}
 }
\caption{Dual norms of the residuals and norms of the global errors with respect to time for different choices of the POD tolerance $tol$. Femoropopliteal bypass application in which the high fidelity solutions are obtained using a finite element approximation on a fine mesh.}
\label{fig:6.2.9}
\end{figure}

{\color{black} To explain why we obtain better results for the pressure than for the velocity, we investigate the absolute values of velocity, pressure and global solutions errors and we compare them to the dual norms of the residuals (see Figure \ref{fig:6.2.9}). Since the velocity and pressure norms have two different magnitudes ($10-10^2$ for the velocity and $10^3-10^5$ for the pressure), the corresponding absolute values of the pressure errors are bigger than the velocity ones, even if the relative errors are lower. The greedy procedure selects the worst approximated time instant based on the dual norms of the residuals and these quantities are indicators of the global absolute errors. Since the latter is mostly due to the pressure error, this can explain why the greedy enrichment provides better results for the pressure than for the velocity.

\section{Conclusions} 
In this work we presented an application of reduced order modelling to a RFSI problem that is indeed an unsteady Navier-Stokes problem with generalized Robin boundary conditions. We detailed how an affine decomposition with respect to boundary data varying in time can be obtained under suitable hypothesis. Moreover, we presented and detailed how the POD can be applied to the RFSI problem in order to take into account the different order of magnitudes of the variables. We discussed the introduction of the supremizer functions inside the reduced basis, necessary to include the pressure in the reduced system. Afterwards, we proposed an enrichment of the POD reduced basis based on a greedy algorithm. All the algorithms presented were then numerically tested on a realistic haemodynamics problem. We tested the POD and greedy enrichment algorithm on two cases: a test case, where the finite element solution is obtained with a coarse grid, and a fine case, where the finite element space has order of $10^6$ degrees of freedom. The results showed the good performances of the POD reduction algorithm on the RFSI problem, also with respect to the evaluation of specific haemodynamics target output (wall shear stress). Moreover we provided numerical evidence of how the reduced approximation can be improved using the greedy enrichment algorithm, in particular regarding the pressure error. The different behaviour of the velocity and pressure errors is due to the use of the dual norm of the residual as an indicator of the global solution error. Indeed, since we do not have suitable a-posteriori error estimators, one for the velocity and one for the pressure variables, we measure the dual norm of the residual as a surrogate estimator. Being the pressure variable and the correspondent error two order of magnitudes grater that the velocity ones, the residual is indeed an indicator of the pressure errors. Nevertheless, even in lack of theoretical results, numerical experiments showed that the greedy enrichment is able to improve the quality of the reduced approximation allowing us to save computational time. The development of suitable a-posteriori error estimators for the pressure and velocity in the case of RFSI problem would be required to improve the performances of the greedy enrichment.

\section*{Acknowledgments}
The authors acknowledge support of the European Research Council under the Advanced Grant ERC-2008-AdG 227058, Mathcard, Mathematical Modelling and Simulation of the Cardiovascular System. Moreover, we thank Prof. Bernard Haasdonk and Prof. Alfio Quarteroni for the insightful discussion. 
We gratefully acknowledge the Swiss National Supercomputing Center (CSCS) for providing the CPU resources for the numerical simulations under the projects IDs s475 and s796. 
% bibliography
\bibliographystyle{siam}
\bibliography{bibliography}

\begin{thebibliography}{10}

\bibitem{Amsallem:2011aa}
{\sc D.~Amsallem and C.~Farhat}, {\em An online method for interpolating linear
  parametric reduced-order models}, SIAM J. Sci. Comput., 33 (2011),
  pp.~2169--2198.

\bibitem{Bazilevs:2007aa}
{\sc Y.~Bazilevs, C.~Michler, V.~M. Calo, and T.~J.~R. Hughes}, {\em Weak
  dirichlet boundary conditions for wall-bounded turbulent flows}, Computer
  Methods in Applied Mechanics and Engineering, 196 (2007), pp.~4853 -- 4862.

\bibitem{Burkardt:2006aa}
{\sc J.~Burkardt, M.~Gunzburger, and H.-C. Lee}, {\em \{POD\} and cvt-based
  reduced-order modeling of {N}avier–{S}tokes flows}, Computer Methods in
  Applied Mechanics and Engineering, 196 (2006), pp.~337 -- 355.

\bibitem{Colciago:2012aa}
{\sc C.~M. Colciago, S.~Deparis, and A.~Quarteroni}, {\em Comparisons between
  reduced order models and full 3{D} models for fluid-structure interaction
  problems in haemodynamics}, J. Comput. Appl. Math., 265 (2014), pp.~120 --
  138.
\newblock Current Trends and Progresses in Scientific Computation Dedicated to
  Professor Ben-yu Guo on His 70th Birthday.

\bibitem{Deparis:2008aa}
{\sc S.~Deparis}, {\em Reduced basis error bound computation of
  parameter-dependent {N}avier-{S}tokes equations by the natural norm
  approach}, SIAM Journal on Numerical Analisys, 46 (2008), pp.~2039--2067.

\bibitem{Deparis:2012ab}
{\sc S.~Deparis and A.~E. Lovgren}, {\em Stabilized reduced basis approximation
  of incompressible three-dimensional {N}avier-{S}tokes equations in
  parametrized deformed domains}, Journal of Scientific Computing, 50 (2012),
  pp.~198--212.

\bibitem{Deparis:2009aa}
{\sc S.~Deparis and G.~Rozza}, {\em Reduced basis method for
  multi-parameter-dependent steady {N}avier–{S}tokes equations: Applications
  to natural convection in a cavity}, Journal of Computational Physics, 228
  (2009), pp.~4359 -- 4378.

\bibitem{Figueroa:2009aa}
{\sc C.~A. Figueroa, S.~Baek, C.~A. Taylor, and J.~D. Humphrey}, {\em A
  computational framework for fluid-solid-growth modeling in cardiovascular
  simulations}, Computer Methods in Applied Mechanics and Engineering, 198
  (2009), pp.~3583--3601.

\bibitem{Figueroa:2006aa}
{\sc C.~A. Figueroa, I.~E. Vignon-Clementel, K.~E. Jansen, T.~J.R. Hughes, and
  C.~A. Taylor}, {\em A coupled momentum method for modeling blood flow in
  three-dimensional deformable arteries}, Computer Methods in Applied Mechanics
  and Engineering, 195 (2006), pp.~5685--5706.

\bibitem{Formaggia:2009aa}
{\sc L.~Formaggia, A.~Quarteroni, and A.~Veneziani}, {\em Cardiovascular
  mathematics}, Springer, 2009.

\bibitem{Gerner:2012aa}
{\sc A.-L. Gerner, A.~Reusken, and K.~Veroy}, {\em Reduced basis a posteriori
  error bounds for the instationary {S}tokes equations}.
\newblock submitted, 2012.

\bibitem{Gerner:2011aa}
{\sc A.-L. Gerner and K.~Veroy}, {\em Reduced basis a posteriori error bounds
  for the {S}tokes equations in parametrized domains: A penalty approach},
  Mathematical Models and Methods in Applied Sciences, 21 (2011),
  pp.~2103--2134.

\bibitem{Giordana:2005aa}
{\sc S.~Giordana, S.~J. Sherwin, J.~Peiró, D.~J. Doorly, J.~S. Crane, K.~E.
  Lee, N.~J. Cheshire, and C.~G. Caro}, {\em Local and global geometric
  influence on steady flow in distal anastomoses of peripheral bypass grafts},
  J Biomech Eng., 127 (2005), pp.~1087--1098.

\bibitem{Grinberg:2009aa}
{\sc L.~Grinberg, A.~Yakhot, and G.~Karniadakis}, {\em Analyzing transient
  turbulence in a stenosed carotid artery by proper orthogonal decomposition},
  Annals of Biomedical Engineering, 37 (2009), pp.~2200--2217.

\bibitem{Gunzburger:2007aa}
{\sc M.~D. Gunzburger, J.~S. Peterson, and J.~N. Shadid}, {\em Reduced-order
  modeling of time-dependent \{PDEs\} with multiple parameters in the boundary
  data}, Computer Methods in Applied Mechanics and Engineering, 196 (2007),
  pp.~1030 -- 1047.

\bibitem{Haasdonk:2008aa}
{\sc B.~Haasdonk and M.~Ohlberger}, {\em Reduced basis method for finite volume
  approximations of parametrized linear evolution equations}, ESAIM:
  Mathematical Modelling and Numerical Analysis, 42 (2008), pp.~277--302.

\bibitem{hesthaven2016certified}
{\sc J.~S. Hesthaven, G.~Rozza, and B.~Stamm}, {\em Certified reduced basis
  methods for parametrized partial differential equations}, Springer, 2016.

\bibitem{Iapichino:2012aa}
{\sc L.~Iapichino, A.~Quarteroni, and G.~Rozza}, {\em A reduced basis hybrid
  method for the coupling of parametrized domains represented by fluidic
  networks}, Computer Methods in Applied Mechanics and Engineering, 221–222
  (2012), pp.~63 -- 82.

\bibitem{Takahito:2014aa}
{\sc T.~Kashiwabara, C.~Colciago, L.~Ded\'e, and A.~Quarteroni}, {\em
  Well-posedness, regularity, and convergence analysis of the finite element
  approximation of a generalized robin boundary value problem}.
\newblock submitted, 2014.

\bibitem{Lassila:2013aa}
{\sc T.~Lassila, A.~Manzoni, A.~Quarteroni, and G.~Rozza}, {\em Boundary
  control and shape optimization for the robust design of bypass anastomoses
  under uncertainty}, ESAIM: Mathematical Modelling and Numerical Analysis, 47
  (2013), pp.~1107--1131.

\bibitem{Lassila:2014aa}
{\sc T.~Lassila, A.~Manzoni, A.~Quarteroni, and G.~Rozza}, {\em Model order
  reduction in fluid dynamics: challenges and perspectives}, in Reduced Order
  Methods for modeling and computational reduction, A.~Quarteroni and G.~Rozza,
  eds., vol.~9 of MS\&A Series, Springer, 2014.

\bibitem{Loth:2008aa}
{\sc F.~Loth, P.~F. Fischer, and H.~S. Bassiouny}, {\em Blood flow in
  end-to-side anastomoses}, Annual Review of Fluid Mechanics, 40 (2008),
  pp.~367--393.

\bibitem{Loth:2002aa}
{\sc F.~Loth, S.~A. Jones, C.~K. Zarins, D.~P. Giddens, R.~F. Nassar,
  S.~Glagov, and H.~S. Bassiouny}, {\em Relative contribution of wall shear
  stress and injury in experimental intimal thickening at {P}{T}{F}{E}
  end-to-side arterial anastomoses}, J Biomech Eng, 124 (2002), pp.~44--51.

\bibitem{Manzoni:2014aa}
{\sc A.~Manzoni and F.~Negri}, {\em {Rigorous and heuristic strategies for the
  approximation of stability factors in nonlinear parametrized PDEs}}.
\newblock Submitted. Preprint available as Mathicse Report Nr. 08.2014., 2014.

\bibitem{Marchandise:2011aa}
{\sc E.~Marchandise, P.~Crosetto, C.~Geuzaine, J.-F. Remacle, and E.~Sauvage},
  {\em Quality open source mesh generation for cardiovascula flow simulation},
  in Modeling of Physiological Flows, D. Ambrosi and A. Quarteroni and G.
  Rozza, 2011.

\bibitem{Moireau:2011aa}
{\sc P.~Moireau, N.~Xiao, M.~Astorino, C.~A. Figueroa, D.Chapelle, C.~A.
  Taylor, and J.F. Gerbeau}, {\em External tissue support and fluid-structure
  simulation in blood flows}, Biomechanics and Modeling in Mechanobiology, 11
  (2012), pp.~1--18.

\bibitem{Nguyen:2009aa}
{\sc N.-C. Nguyen, G.~Rozza, and A.T. Patera}, {\em Reduced basis approximation
  and a posteriori error estimation for the time-dependent viscous burgers’
  equation}, Calcolo, 46 (2009), pp.~157--185.

\bibitem{Nobile:2008aa}
{\sc F.~Nobile and C.~Vergara}, {\em An effective fluid-structure interaction
  formulation for vascular dynamics by generalized robin conditions}, SIAM
  Journal on Scientific Computing, 30 (2008), pp.~731--763.

\bibitem{quarteroni2015reduced}
{\sc A.~Quarteroni, A.~Manzoni, and F.~Negri}, {\em Reduced basis methods for
  partial differential equations: an introduction}, vol.~92, Springer, 2015.

\bibitem{CECAM:2013aa}
{\sc A.~Quarteroni and G.~Rozza}, eds., {\em Reduced {O}rder {M}ethods for
  {M}odeling and {C}omputational {R}eduction}, MS\&A, Modeling, Simulation and
  Applications, Springer, Milano, 2013.
\newblock Book edited with selected contributions from CECAM workshop, EPFL,
  May 2012.

\bibitem{Rowley:2005aa}
{\sc C.~W. Rowley}, {\em Model reduction for fluids, using balanced proper
  orthogonal decomposition}, International Journal of Bifurcation and Chaos
  (IJBC), 15 (2005), pp.~997--1013.

\bibitem{Rozza:2008aa}
{\sc G.~Rozza, D.~B.~P. Huynh, and A.~T. Patera}, {\em Reduced basis
  approximation and a posteriori error estimation for affinely parametrized
  elliptic coercive partial differential equations}, Archives of Computational
  Methods in Engineering, 15 (2008), pp.~229--275.

\bibitem{Rozza:2007aa}
{\sc G.~Rozza and K.~Veroy}, {\em On the stability of the reduced basis method
  for {S}tokes equations in parametrized domains}, Computer Methods in Applied
  Mechanics and Engineering, 196 (2007), pp.~1244 -- 1260.

\bibitem{Sen:2006aa}
{\sc S.~Sen, K.~Veroy, D.~B.~P. Huynh, S.~Deparis, N.~C. Nguyen, and A.~T.
  Patera}, {\em “{N}atural norm” a posteriori error estimators for reduced
  basis approximations}, Journal of Computational Physics, 217 (2006), pp.~37
  -- 62.

\bibitem{Sirisup:2005aa}
{\sc S.~Sirisup and G.~E. Karniadakis}, {\em Stability and accuracy of periodic
  flow solutions obtained by a pod-penalty method}, Physica D: Nonlinear
  Phenomena, 202 (2005), pp.~218 -- 237.

\bibitem{Willcox:2002aa}
{\sc K.~Willcox and J.~Peraire}, {\em Balanced model reduction via the proper
  orthogonal decomposition}, AIAA Journal,  (2002), pp.~2323--2330.

\bibitem{Wirtz:2014aa}
{\sc D.~Wirtz, D.C. Sorensen, and B.~Haasdonk}, {\em A-posteriori error
  estimation for deim reduced nonlinear dynamical systems}, SIAM J. Sci. Comp.,
  36 (2014), pp.~A311--A338.

\bibitem{Yano:2012aa}
{\sc Masayuki Yano, Anthony~T. Patera, and Karsten Urban}, {\em A space-time
  hp-interpolation-based certified reduced basis method for {B}urgers'
  equation}, Mathematical Models and Methods in Applied Sciences, 0 (2014),
  pp.~1--33.

\end{thebibliography}
% -------------------------------------------------------- EOD --------------------------------------------------------
\end{document}